\documentclass[preprint,12pt,times]{elsarticle}
\usepackage{lineno}
\modulolinenumbers[200]
\usepackage{color}
\usepackage{colortbl}
\usepackage{transparent}
\usepackage{graphicx}
\usepackage{import}
\usepackage{diagbox}
\usepackage{algorithmic}
\usepackage{algorithm}
\usepackage{appendix}
\usepackage{amsmath}
\usepackage{amscd}
\usepackage{mathrsfs}
\usepackage{graphicx,color}
\usepackage{arydshln}
\usepackage{booktabs}
\usepackage{multirow}
\usepackage{rotating}
\usepackage{setspace}
\usepackage{indentfirst}
\usepackage{cases}
\usepackage{verbatim}
\usepackage{xcolor}
\graphicspath{
	{./figures/}
	{./figures/2d_forward/},
	{./figures/3d_forward/},
	{./figures/3d_IP/}
}
\usepackage{amssymb,color}
\usepackage{amsthm}
\usepackage{epsfig,subfigure}

\newcommand{\edge}[1]{\ar@{-}[#1]}
\numberwithin{equation}{section}

\biboptions{numbers,sort&compress}%

\begin{document}

\begin{frontmatter}
    \title{Laplace-fPINNs: Laplace-based fractional physics-informed neural networks for solving forward and inverse problems of subdiffusion}

    \author[1,4]{Xiong-Bin Yan}
    \ead{yanxb2015@163.com}
    \author[1,2,3]{Zhi-Qin John Xu\corref{cor}}
    \ead{xuzhiqin@sjtu.edu.cn}
    \author[1,2,3,4]{Zheng Ma\corref{cor}}
    \ead{zhengma@sjtu.edu.cn}
    \cortext[cor]{Corresponding author}
    \address[1]{School of Mathematical Sciences, Shanghai Jiao Tong University, Shanghai, China}
    \address[2]{Institute of Natural Sciences, MOE-LSC, Shanghai Jiao Tong University, Shanghai, China}
    \address[3]{Qing Yuan Research Institute, Shanghai Jiao Tong University, Shanghai, China}
    \address[4]{CMA-Shanghai, Shanghai Jiao Tong University, Shanghai, China}
    \begin{abstract}
        The use of Physics-informed neural networks (PINNs) has shown promise
        in solving forward and inverse problems of fractional diffusion equations.
        However, due to the fact that automatic differentiation is not applicable for
        fractional derivatives, solving fractional diffusion equations using PINNs
        requires addressing additional challenges. To address this issue, this paper proposes an extension to PINNs called
        Laplace-based fractional physics-informed neural networks (Laplace-fPINNs),
        which can effectively solve the forward and inverse problems of fractional
        diffusion equations. This approach avoids introducing a mass of auxiliary
        points and simplifies the loss function. We validate the effectiveness of the Laplace-fPINNs approach using several examples. Our numerical results demonstrate that the Laplace-fPINNs method can effectively solve both the forward and inverse problems of high-dimensional fractional diffusion equations.
    \end{abstract}
\end{frontmatter}

\section{Introduction}\label{Sec1}
Fractional diffusion equations have been studied extensively in engineering, physics, and mathematical literature owing to their superior capability for modeling anomalous
diffusion phenomena. The model differs from the standard diffusion model in that it follows
a basic assumption that the diffusion obeys the standard Brownian motion and has been applied
to animal coat patterns and nerve cell signals. A distinctive feature of standard Brownian motion is that the mean squared displacement $\langle x^2(t)\rangle$ of diffusing species increases linearly with time, i.e., $\langle x^2(t)\rangle\sim K_1t$. However, in anomalous diffusion, the mean squared displacement shows a non-linear power law growth with time, i.e., $\langle x^2(t)\rangle\sim K_{\alpha}t^{\alpha}$, where $0<\alpha<1$ represents subdiffusion, and $\alpha>1$ represents superdiffusion.  At the microscopic level, such anomalous diffusion processes can be accurately described by a continuous-time random walk where the waiting time between successive particle leaps follows a heavy-tailed distribution with a diverging mean. At the macroscopic level, anomalous diffusion describes the evolution of the probability density function of a particle that appears at a given spatial location $x$ and time $t$. A list of successful applications of fractional diffusion equations are extensive and continually expanding. These applications include, but are not limited to, solute transport in heterogeneous media~\cite{Berkowitz_Cortis_Dentz_Scher_2006, Dentz_Cortis_Scher_Berkowitz_2004}, thermal diffusion on fractal domains~\cite{Nigmatullin_1986}, protein transport within membranes~\cite{Kou_2008, Ritchie_Shan_Kondo_Iwasawa_Fujiwara_Kusumi_2005}, and flow in highly heterogeneous aquifers~\cite{Berkowitz_Klafter_Metzler_Scher_2002}. Comprehensive reviews of the physics modeling and a diverse range of applications can be found in~\cite{Bouchaud+Georges-1990, Metzler+Klafter-2000}.

In this paper, we consider the following time-fractional diffusion equation on a bounded domain $\Omega\subset \mathbb{R}^{d}$ with homogeneous Dirichlet boundary condition
\begin{align}\label{1.1}
    \left \{
    \begin{aligned}
         & \partial_{0+}^{\alpha}u(x,t)= \nabla\cdot(a(x) \nabla u(x,t)) + c(x)u(x,t) + f(x,t),
        ~ x\in\Omega,~t>0,                                                                      \\
         & u(x,0)=u_0(x),~x\in \Omega,                                                          \\
         & u(x,t)=0,~x\in\partial \Omega,~ t>0.
    \end{aligned}
    \right.
\end{align}
The left-hand side of the above equation is a Caputo derivative of order $\alpha$, which is defined by
\begin{align*}
    \partial_{0+}^{\alpha}u(x,t)=\frac{1}{\Gamma(1-\alpha)}\int_{0}^{t}{(t-\tau)}^{-\alpha}\frac{\partial u(x,\tau)}{\partial \tau}d\tau, ~0<\alpha<1,
\end{align*}
where $\Gamma(\cdot)$ is the gamma function, $a(x)$, $c(x)$, $f(x,t)$ and $u_0(x)$ represent the diffusion coefficient, reaction coefficient, source and initial value,
respectively.

Solving partial differential equations numerically is a well-known challenge, and it becomes even more difficult when dealing with fractional diffusion equations that involve nonlocal operators. Recently, there has been a growing trend to apply machine learning techniques for solving forward and inverse problems of partial differential equations. Several examples include the use of Gaussian process regression~\cite{Graepel-2003,Raissi+Perdikaris+Karniadakis-2018,Lee+Bilionis+Tepole-2020} and deep learning-based methods~\cite{Lagaris+Likas+Fotiadis-1998,Raissi+Perdikaris+Karniadakis-2019,Lu+Meng+Mao+Karniadakis-2021,E+Yu-2018,Zang_Bao_Ye_Zhou_2020,Bao_Ye_Zang_Zhou_2020} to solve these types of problems.

This paper focuses on the use of deep learning methods for solving partial differential equations (PDEs), which can be categorized into two approaches. The first approach utilizes neural networks to learn the solution operator for a given problem using a large number of input-output pairs obtained from numerical simulation or governing equations. In this case, deep neural networks are used to approximate a function that maps the input (e.g.\ coefficients, initial value, source, and boundary) to the solution of partial differential equations. For more detailed research, see~\cite{Lu_Jin_Karniadakis_2019,Jin_Meng_Lu_2022,zhu_Zabaras_2018,Li_Kovachki_Azizzadenesheli_liu_Bhattacharya_Stuart_Anandkumar_2021,Li_Kovachki_Azizzadenesheli_Liu_Bhattacharya_Anandkumar_2020,Li_Kovachki_Azizzadenesheli_Liu_Bhattacharya_Stuart_Anandkumar_2020_1,zhang2022mod,yan2022bayesian}. The second approach utilizes deep neural networks to approximate the solution of PDEs due to their high expressive power. This method takes advantage of recent advances in automatic differentiation, one of the useful techniques in scientific computing, to derive the derivative of neural networks to obtain the cost function, the trainable parameters of the neural network can be optimized by minimizing the cost function. This approach is simple yet powerful, and introduces potentially physics-informed deep learning methods for solving PDEs. Examples of this approach include Deep Ritz methods~\cite{E+Yu-2018, Liao+Ming-2021}, Physics-informed neural networks methods~\cite{Lu+Meng+Mao+Karniadakis-2021, Raissi+Perdikaris+Karniadakis-2019, Lou_Meng_Karniadakis_2021, Kharazmi_Zhang_Karniadakis_2021, Yang_Meng_Karniadakis_2021,Shukla+Jagtap+Karniadakis-2021,Mao+Lu+Marxen+Zaki+Karniadakis-2021}, weak adversarial networks methods~\cite{Zang_Bao_Ye_Zhou_2020, Bao_Ye_Zang_Zhou_2020}, Deep Galerkin methods~\cite{Sirignano+Spiliopoulos-2018, Shang+Wang+Sun-2022}, and multi-scale DNN~\cite{liu2020multi,li2020elliptic}.

Recently, solving fractional partial differential equations by neural networks attracts more and more attention~\cite{Pang+Lu+Karniadakis-2019,Guo+Wu+Yu+Zhou-2022,yan2022bayesian}. However, the classical chain rule is not applicable in fractional calculus, which renders automatic differentiation ineffective for fractional derivatives. To overcome this challenge, Pang et al.\ \cite{Pang+Lu+Karniadakis-2019} propose an extension of Physics-informed neural networks (PINNs) called fractional physics-informed neural networks, which approximates fractional derivatives by numerically discretizing fractional operators and employs automatic differentiation to analytically obtain the integer-order derivatives of neural networks. Similarly, Guo et al.\ \cite{Guo+Wu+Yu+Zhou-2022} propose a Monte Carlo sampling-based PINN method that computes the fractional derivatives of neural network output via Monte Carlo sampling. However, both methods require the introduction of auxiliary points to calculate the fractional derivative, which increases the computational cost of neural network training.

In this paper, we propose a Laplace-based fractional physics-informed neural network, named Laplace-fPINNs, for solving forward and inverse problems of subdiffusion. The proposed method involves transforming the original time-fractional diffusion equation into a restricted equation in Laplace space using the Laplace transform, which can then be solved using physics-informed neural networks. The restricted equation in Laplace space does not contain fractional derivatives, which circumvents the issue of applying automatic differentiation to the time-fractional derivative.  A numerical inverse Laplace transform to convert the PINNs solution from Laplace space to the time domain. In summary, the
main contributions of this paper are listed as follows:

\begin{itemize}
    \item We provide a novel approach, called Laplace-fPINNs, for solving the time-fractional diffusion equation (\ref{1.1}). This method substantially reduces the computational cost of the loss function during neural network training and eliminates the requirement for auxiliary points used in previous studies~\cite{Pang+Lu+Karniadakis-2019,Guo+Wu+Yu+Zhou-2022} to approximate the time-fractional derivative.

    \item We utilize the proposed Laplace-fPINNs method to tackle the challenging task of identifying a high-dimensional diffusion coefficient from given measurements. This problem is typically difficult to solve using traditional inversion methods.
\end{itemize}

The structure of this paper is as follows. In Section~\ref{sec02}, we review the fundamentals of physics-informed neural networks (PINNs), and summarize some previous research on using PINNs to solve fractional partial differential equations. In Section~\ref{sec2}, we elaborate
Laplace-based physics-informed neural networks for solving the time-fractional diffusion equation (\ref{1.1}). Section~\ref{sec3} presents a comprehensive numerical investigation of the performance of the Laplace-fPINNs approach for solving both forward and inverse problems associated with equation (\ref{1.1}). Finally, we conclude
the paper in Section~\ref{sec4}.

\section{Prelimilaries}\label{sec02}
\subsection{PINNs}
\begin{figure}[H]
    \centering
    \includegraphics[height=3.5cm,width=11cm]{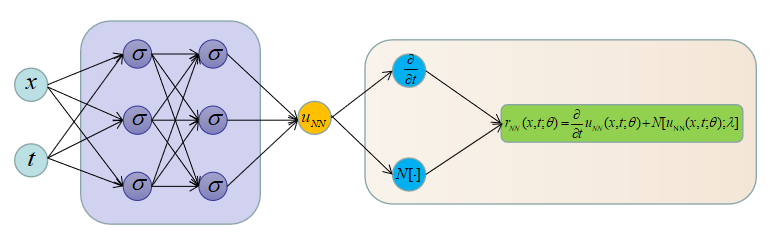}
    \caption{PINNs for solving partial differential equations.}\label{pinns}
\end{figure}

In this part, we first recall the main idea of physics-informed neural networks (PINNs) for solving the integer-order nonlinear partial differential equations of general form
\begin{align}\label{2.1}
    u_{t}+\mathcal{N}[u;\lambda]=0, ~x\in \Omega,~t\in [0,T],
\end{align}
where $u(x,t)$ denotes the solution of equation (\ref{2.1}), $\mathcal{N}[\cdot;\lambda]$ is a nonlinear operator
parameterized by $\lambda$, $\Omega$ is a bounded domain in $\mathbb{R}^{d}$.
The PINNs approach solves the above equation (\ref{2.1}) by approximating $u(x,t)$ with a deep neural network, which takes the coordinate $(x,t)$ as input and outputs a scalar $u_{NN}(x,t;\theta)$, where $\theta$ represents all trainable parameters of the neural network. The approximate solution $u_{NN}(x,t;\theta)$ is differentiable and can be substituted into equation (\ref{2.1}) using the automatic differentiation~\cite{Baydin+Pearlmutter+Radul+Siskind-2017} to obtain
\begin{align}\label{2.2}
    r_{NN}(x,t;\theta):=\partial_t u_{NN}(x,t;\theta)+\mathcal{N}[u_{NN}(x,t;\theta);\lambda].
\end{align}
Both neural networks $u_{NN}(x,t;\theta)$ and $r_{NN}(x,t;\theta)$ share the same trainable parameters which can be learned by minimizing the following loss
\begin{align*}
    L(\theta) & := L_{bd}(\theta)+L_{eq}(\theta)                                             \\
              & = \frac{1}{N_u}\sum_{i=1}^{N_u}|u_{NN}(x_{u}^{i},t_{u}^{i};\theta)-u^{i}|^2+
    \frac{1}{N_{r}}\sum_{i=1}^{N_r}|r_{NN}(x_{r}^{i}, t_{r}^{i};\theta)|^2,
\end{align*}
where $\{x_{u}^{i},t_{u}^{i},u^{i}{\}}_{i=1}^{N_u}$ represent the initial and boundary training data on $u(x,t)$, $\{x_{r}^{i}, t_{r}^{i}{\}}_{i=1}^{N_r}$ specify the collocation points for $r_{NN}(x,t;\theta)$. Figure~\ref{pinns} shows the sketch of the PINNs.

\subsection{Using PINNs to solve fractional PDEs}
PINNs are a new approach to solve partial differential equations (PDEs) using neural networks. However, when dealing with fractional PDEs, which involve non-integer order derivatives, PINNs face some challenges. One of the primary challenges is the computation of the fractional Caputo derivative, which is widely used in physics and engineering. Unlike integer order derivatives, Caputo derivatives depend on the entire history of the function, which means that the numerical integration becomes more expensive and memory-intensive as the simulation time progresses. Moreover, Caputo derivatives do not obey the classical chain rule, which implies that the standard technique of auto-differentiation in PINNs cannot be directly applied. Therefore, new methods and algorithms are needed to overcome these difficulties and make PINNs more efficient and robust for fractional PDEs.

Inspired by the previous PINNs approach for solving the integer-order partial differential equations, Pang et al.\ in~\cite{Pang+Lu+Karniadakis-2019} propose a fractional PINNs (fPINNs) approach for solving fractional partial differential equations. Specifically, fPINNs employ automatic differentiation to analytically derive the integer-order derivatives of neural network outputs and use numerical discretization to approximate the fractional derivatives. Moreover, to reduce the excessive costs of the numerical discretization for the fPINNs, Guo et al.\ in~\cite{Guo+Wu+Yu+Zhou-2022} propose a Monte Carlo sampling-based PINNs (MC-fPINNs) approach for solving fractional partial differential equations. This method simplifies the time-fractional derivative as an expectation of a function, which is solved by the Monte Carlo method by inducing a probability density function. Although both fPINNs and MC-fPINNs have demonstrated high accuracy in solving forward and inverse problems of fractional partial differential equations, they require auxiliary points to derive the time-fractional derivative of the neural network outputs, increasing the computational cost of approximating fractional derivatives by neural networks.

\section{Laplace-fPINNs}\label{sec2}
To tackle the above problems, we propose a novel method called \textit{Laplace-based fractional physics-informed neural networks (L-fPINNs)} that leverages the Laplace transform to compute the time-fractional derivative of the neural network outputs analytically. This method avoids the need for auxiliary points and numerical discretization, and reduces the computational cost and complexity of solving fractional partial differential equations by neural networks. The proposed method consists of three ingredients: firstly, the Laplace transform of the fractional PDE;\@ secondly, the PINNs method in Laplace domain and finally, the numerical inversion to get the final solution. We show that by applying the Laplace transform, we can reduce the fractional PDE to an ordinary PDE in Laplace domain, which can be solved efficiently by PINNs. Then, we use a numerical inversion technique to recover the solution in the original domain.

\subsection{Laplace transform}
First, we present a definition of the one-dimensional Laplace transform which converts a function from a time domain to a Laplace space domain. Let $f(t)$ be a dependent on time variable $t\ge 0$ function, $\tilde{f}(s)$ represents its image in Laplace space, which defined by
\begin{align}\label{2.3}
    \tilde{f}(s)=\mathcal{L}[f(t)](s)=\int_{0}^{\infty}e^{-st}f(t)dt,\quad (s\in\mathbb{C}).
\end{align}
If the integral (\ref{2.3}) is convergent at a point $s_0\in\mathbb{C}$, then it converges absolutely
for any $s\in\mathbb{C}$ such that $\Re(s)>\Re(s_0)$, where $\Re(s)$ represents the real part of complex
variable s. The Laplace transform of time-fractional Caputo derivative $\partial_{0+}^{\alpha}f(t)$ can be written as
\begin{align*}
    \mathcal{L}[\partial_{0+}^{\alpha}f](s)
     & =s^{\alpha}\mathcal{L}[f(t)](s)-s^{\alpha-1}f(0)                   \\
     & =s^{\alpha}\tilde{f}(s)-s^{\alpha-1}f(0),\quad t\ge 0,~0<\alpha<1.
\end{align*}

Based on the above introduction, we apply the Laplace transform on both sides of equation (\ref{1.1}), leading to
\begin{align}\label{2.4}
    \left \{
    \begin{aligned}
         & s^{\alpha}\tilde{u}(x,s) - \nabla\cdot(a(x) \nabla \tilde{u}(x,s)) - c(x)\tilde{u}(x,s) = s^{\alpha-1}u_0(x) +\tilde{f}(x,s),~x\in\Omega,~s\in\mathbb{C}, \\
         & \tilde{u}(x,s)=0,~x\in\partial \Omega,~s\in\mathbb{C},
    \end{aligned}
    \right.
\end{align}
where $\tilde{f}(x,s)=\mathcal{L}[f(x,t)](x,s)=\int_{0}^{\infty}e^{-st}f(x,t)dt$. It is worth noting that the time-fractional diffusion equation in Laplace space no longer contains the time-fractional derivative $\partial_{0+}^{\alpha}u(x,t)$,  as opposed to the original fractional diffusion equation (\ref{1.1}). This implies that if we solve the numerical solution $u$ directly via equation (\ref{2.4}), we can substantially reduce the computational expense associated with discretization of the time-fractional derivative.

\subsection{PINNs in Laplace domain}
Next, we propose a Lapalce-based physics-informed neural networks approach to solve equation (\ref{2.4}) approximately. We reformulate equation (\ref{2.4}) as follows:
\begin{align}\label{2.5}
    \left \{
    \begin{aligned}
         & \mathcal{A}[\tilde{u}](x,s) = F(x,s), ~x\in\Omega,~s\in\mathbb{C}, \\
         & \tilde{u}(x,s)=0,~x\in\partial \Omega,~s\in\mathbb{C},
    \end{aligned}
    \right.
\end{align}
where $\mathcal{A}[\tilde{u}]=s^{\alpha}\tilde{u}(x,s)-\nabla\cdot(a(x)\nabla \tilde{u}(x,s))-c(x)\tilde{u}(x,s)$,
$F(x,s)=s^{\alpha-1}u_0(x)+\tilde{f}(x,s)$. We employ a deep neural network $\tilde{u}_{NN}(x,s;\theta)$ to approximate the solution $\tilde{u}$ in equation (\ref{2.5}), where $\theta$ represents all trainable parameters of the neural network. Subsequently, we substitute $\tilde{u}_{NN}(x,s;\theta)$ into equation (\ref{2.5}) to obtain the corresponding residual
\begin{align*}
    r_{NN}(x,s;\theta):=\mathcal{A}[\tilde{u}_{NN}](x,s) - F(x,s).
\end{align*}
The Laplace-based physics-informed neural network can be trained by minimizing the following composite
loss function
\begin{align}\label{2.6}
    L^{lp}(\theta):=w_{bd} L_{bd}^{lp}(\theta)+w_{eq} L_{eq}^{lp}(\theta),
\end{align}
where
\begin{align*}
    L_{eq}^{lp}(\theta) & =\frac{1}{N_{r}}\sum_{i=1}^{N_r}|r_{NN}(x_r^{i}, s_{r}^{i};\theta)|^2,             \\
    L_{bd}^{lp}(\theta) & =\frac{1}{N_{bd}}\sum_{i=1}^{N_{bd}}|\tilde{u}_{NN}(x_{bd}^i, s_{bd}^i;\theta)|^2.
\end{align*}
The batch sizes $N_{r}$ and $N_{bd}$ are used to denote the number of collocation points for the residual $r_{NN}(x,s;\theta)$, $x\in\Omega$ and the boundary $\tilde{u}_{NN}(x,s;\theta)$, $x\in\partial\Omega$, respectively. At each iteration of the stochastic descent algorithm, these collocation points are randomly sampled within the computational domain and boundary. Furthermore, the weight coefficients $w_{eq},~w_{bd}$ in the loss function $L^{lp}(\theta)$ can balance the different learning rates of each loss term. The specific value of the weight coefficients $w_{eq},~w_{bd}$ and the batch sizes $N_{r},~N_{bd}$ will be provided in the experiments.

It is noteworthy that the loss function presented in (\ref{2.6}) excludes the time-fractional Caputo derivative $\partial_{0+}^{\alpha}$, and only contains the term $s^{\alpha}\tilde{u}_{NN}(x,s;\theta)-s^{\alpha-1}u_0(x)$. This implies that there is no need to discretize the time-fractional derivative $\partial_{0+}^{\alpha}$ during the training process for minimizing the loss function in (\ref{2.6}). Consequently, there is no requirement for auxiliary points, as suggested in previous studies~\cite{Pang+Lu+Karniadakis-2019,Guo+Wu+Yu+Zhou-2022}.

\subsection{Numerical inverse of Laplace transform}
By minimizing the loss function (\ref{2.6}) using the random gradient descent algorithm, we obtain an approximate solution of equation (\ref{2.5}) denoted by $\tilde{u}_{NN}(x,s;\theta^{*})$. Next, we discuss how to convert the approximate solution $\tilde{u}_{NN}(x,s;\theta^{*})$ in Laplace space domain to the solution $u$ of the original problem (\ref{1.1}), which requires the use of the numerical inverse Laplace transform (NILT). The Laplace transform inversion is a well-known problem that is notoriously ill-conditioned. Numerical inversion can be an unstable process, with difficulties arising from high sensitivity to round-off errors. To address this challenge, researchers have proposed various algorithms for numerically inverting the Laplace transform, such as the Fourier series method \cite{davies1979numerical}, the Talbot algorithm \cite{abate2006unified}, and the Gaver-Stehfest algorithm \cite{Gaver-1966}. In this study, we have opted to use the Gaver-Stehfest NILT algorithm, which is a well-established technique that is easy to implement. Moreover, the summation weights and nodes of the inversion algorithm do not depend on complex numbers. By leveraging this algorithm, we are able to transform the approximate solution $\tilde{u}_{NN}(x,s;\theta^{*})$ from the Laplace space domain to the time-dependent solution $u$ in equation (\ref{1.1}). An approximate solution $u_{NN}(x,t;\theta^{*})$ in the equation (\ref{1.1}) at a specific time $t$ from the Stehfest NILT algorithm is given by
\begin{align}\label{2.7}
    u_{NN}(x,t;\theta^{*})=\frac{\ln 2}{t}\sum_{i=1}^{M}\mu_{i}\tilde{u}_{NN}(x,\frac{\ln 2}{t}i;\theta^{*}),
\end{align}
where the coefficients $\mu_{i }$ given by
\begin{align*}
    \mu_{i}={(-1)}^{\frac{M}{2}+i}\sum_{k=[\frac{i+1}{2}]}^{\min(i,\frac{M}{2})}
    \frac{k^{\frac{M}{2}}(2k)!}{(\frac{M}{2}-k)!(k)!(k-1)!(i-k)!(2k-i)!},
\end{align*}
in which $[C]$ denotes the nearest integers less than or equal to $C$, M in (\ref{2.7}) is a even number.

According to (\ref{2.7}), we can obtain the approximate solution $u_{NN}(x,t;\theta^{*})$ for a given $t$ by taking a linear combination of the sequence ${\{\tilde{u}_{NN}(x,s_i;\theta^{*})\}}_{i=1}^{M}$, where $s_i=\frac{\ln 2}{t}i$ for $i=1,\ldots,M$. The sequence ${\{\tilde{u}_{NN}(x,s_i;\theta^{*})\}}_{i=1}^{M}$ is obtained by inputting the collocation
points ${\{(x,s_i)\}}_{i=1}^{M}$ into the neural network
approximate solution $\tilde{u}_{NN}(x,s;\theta^{*})$, which is a function of variables $x$ and $s$.


It is worth noting that, according to the definition of the Laplace transform (\ref{2.3}), the neural network output $\tilde{u}_{NN}(x,s;\theta^{*})$ needs to approximate the solution $\tilde{u}(x,s)$, $x\in\Omega,~s\in\mathbb{C}$ in equation (\ref{2.5}). However, the choice of training points $s\in\mathbb{C}$ is not necessary for the neural network training process. Actually, to obtain the approximate solution $u_{NN}(x,t;\theta^{*})$ using the Stehfest NILT method in (\ref{2.7}), we only need to evaluate the neural network output $\tilde{u}_{NN}(x,s_i;\theta^{*})$ at $s_i=\frac{\ln 2}{t}i$, for $i=1,\ldots,M$. Therefore, during neural network training, we choose the training points $s\in\mathcal{S}\subset\mathbb{R}^{+}$, which is a bounded closed interval in $\mathbb{R}^{+}$ that will be specified in the following section.

\section{Experiments}\label{sec3}

In this section, we aim to verify the effectiveness of the proposed Laplace-fPINNs method through several numerical examples. Firstly, we evaluate the performance of the Laplace-fPINNs method in solving the forward problem of subdiffusion (\ref{1.1}) in both two-dimensional and three-dimensional cases. Then, we demonstrate the efficiency of the Laplace-fPINNs method in solving an inverse problem of identifying a diffusion coefficient in a three-dimensional subdiffusion problem.
To quantify the accuracy of the Laplace-fPINNs method, we consider a relative $l_2$ error, which is given by
\begin{align*}
    \text{Relative $l_2$ error} = \frac{\|u(x,t)-u_{NN}(x,t)\|_{2}}{\|u(x,t)\|_{2}},
\end{align*}
where $u$ and $u_{NN}$ represent the exact and approximate solutions, respectively.

Throughout all experiments, we use a fully connected neural network with the Swish activation function to approximate the solution $\tilde{u}(x,s)$ of equation (\ref{2.5}). The initial values of the trainable parameters of the neural network are set to default values in all numerical experiments. The neural network is trained using the Adam optimizer with a learning rate of $1\times 10^{-4}$. The hyperparameters of the optimizer are chosen based on default recommendations.

\subsection{2D forward problem}
\subsubsection{2D forward problem with $T=1$}
To study the performance of the Laplace-fPINNs approximation, we consider solving a forward problem of equation (\ref{1.1}) in two-dimensional case. Without loss of generality, we take the domain $\Omega={[0,1]}^2$, $T=1$. The diffusion coefficient and reaction coefficient are assumed to be $a(x,y)=1$ and $c(x,y)=0$, respectively. The initial value and the source are given by
\begin{align*}
    u_0(x,y) & =\sin(\pi x)\sin(\pi y),    \\
    f(x,y)   & =5\sin(2\pi x)\sin(3\pi y).
\end{align*}
We obtain a reference solution $u$ using a $L_1$-type finite difference method~\cite{Sun+Wu-2006}. The grid resolutions of time and space variables in the numerical method are taken as $101$ and $101\times 101$, respectively.

In this experiment, we set the number of residual points and bound points used to compute the equation loss $L_{eq}^{lp}(\theta)$ and the bound loss $L_{bd}^{lp}(\theta)$ for each mini-batch as 1000 and 1600, respectively. For each iteration, the residual points and bound points of the spatial variable are sampled from the uniform distribution on $\Omega$ and $\partial\Omega$, respectively. The training point $s$ is sampled from the interval $\mathcal{S}=[s_{\min},s_{\max}]$, where $s_{\min}=\frac{\ln 2}{T}$, $s_{\max}=\frac{\ln 2}{t_1}M$, $t_1=0.01$. The weights $w_{eq}$ and $w_{bd}$ in (\ref{2.6}) are hyper-parameters that can be manually or automatically tuned during training. The magnitude of $w_{eq}$ and $w_{bd}$ significantly impacts the convergence rate of the neural network training. However, discussing how the weights $w_{eq}$ and $w_{bd}$ affect the convergence rate of the Laplace-fPINNs is beyond the scope of this paper. Thus, in this example, we set the weights $w_{eq}=1$ and $w_{bd}=2000$ heuristically.

We first discuss how the parameter $M$ in (\ref{2.7}), which is used in the Stehfest NILT algorithm, affects the accuracy of the Laplace-fPINNs method. To this end, we fix the number of hidden layers and the number of neurons in each hidden layer of the neural network as 5 and 256, respectively, and set the iteration number to $1.2\times 10^{5}$. To reduce the influence of random initialization of the neural network's parameters on the numerical accuracy, we run the Laplace-fPINNs code five times and report the mean and standard deviation of the relative $l_2$ errors, which are presented in Table~\ref{table_1}.

As shown in Table~\ref{table_1}, our results indicate that the optimal value of $M$ for achieving the best numerical accuracy is 4, and increasing $M$ does not reduce the relative $l_2$ error. This phenomenon can be attributed to the fact that the approximate solution $\tilde{u}_{NN}(x,s;\theta)$ is not accurate enough to the true solution of equation (\ref{2.5}), and the numerical inverse Laplace transform process may introduce additional errors. Similar observations have been reported in prior studies, such as~\cite{Fu+Chen+Yang-2013,Chen+Rashed+Golberg-1998}.

To verify the sensitivity of the Laplace-fPINNs method to the depth and width of the neural network, we calculate the relative $l_2$ error between the approximate solution $u_{NN}(x,t;\theta^{*})$ and the reference solution $u(x,t)$ obtained using the $L_1$-type finite difference method. The value of $M$ is set to the optimal value of 4, as determined in the previous analysis. For each neural network architecture, we run the Laplace-fPINNs code five times. The corresponding numerical results are shown in Table~\ref{table_2}, which includes the mean and standard derivation of the relative $l_2$ error. It can be observed from the table that the width of the neural network has a greater impact on the solution accuracy than the depth. Increasing the width of the neural network results in a slightly higher accuracy.

\begin{figure}[H]
    \centering
    \subfigure[The reference solution $u$.]{
        \begin{minipage}[t]{0.3\linewidth}
            \centering
            \includegraphics[width=2in]{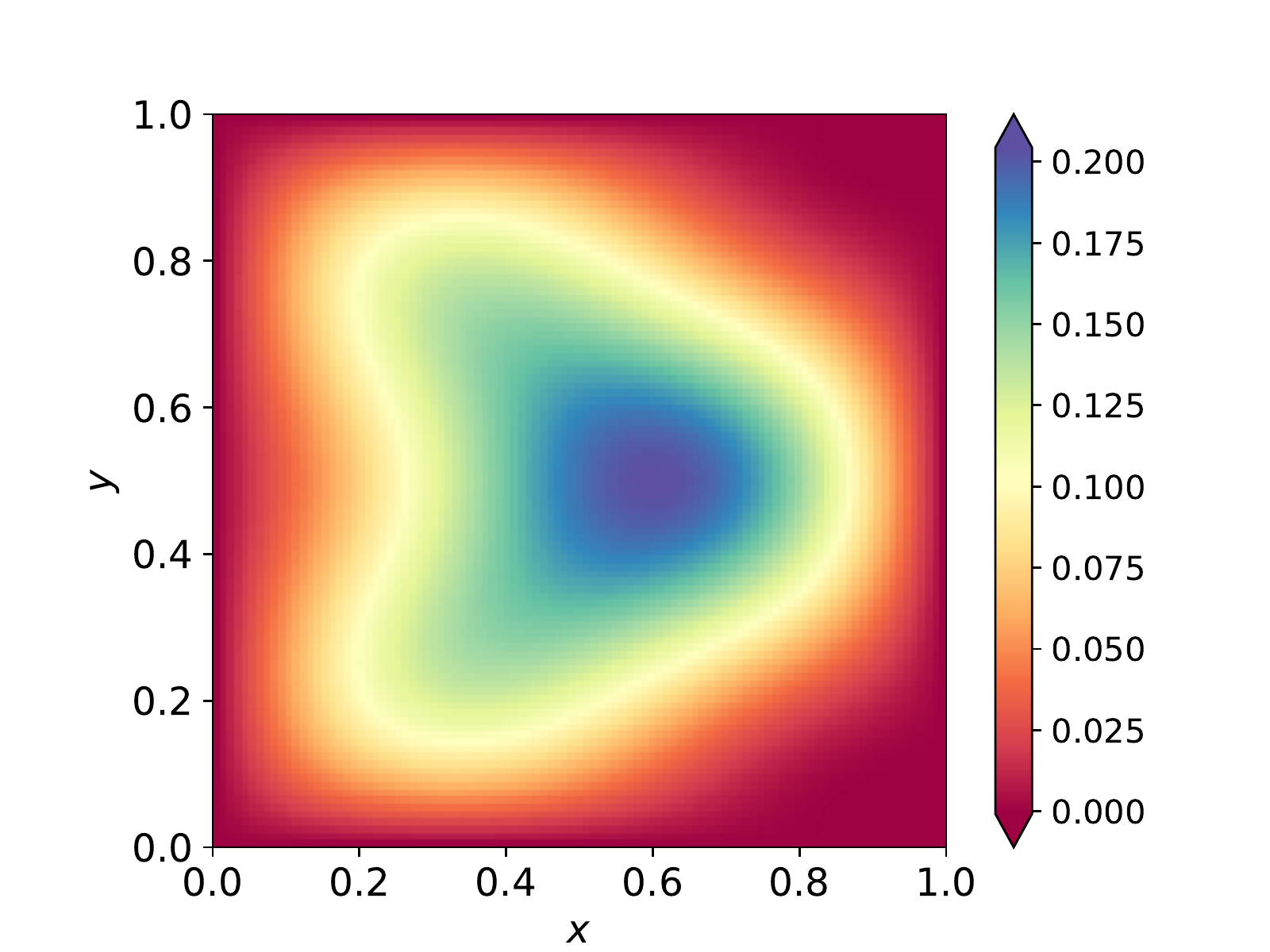}
        \end{minipage}
    }%
    \subfigure[The approximate solution $u_{NN}$.]{
        \begin{minipage}[t]{0.3\linewidth}
            \centering
            \includegraphics[width=2in]{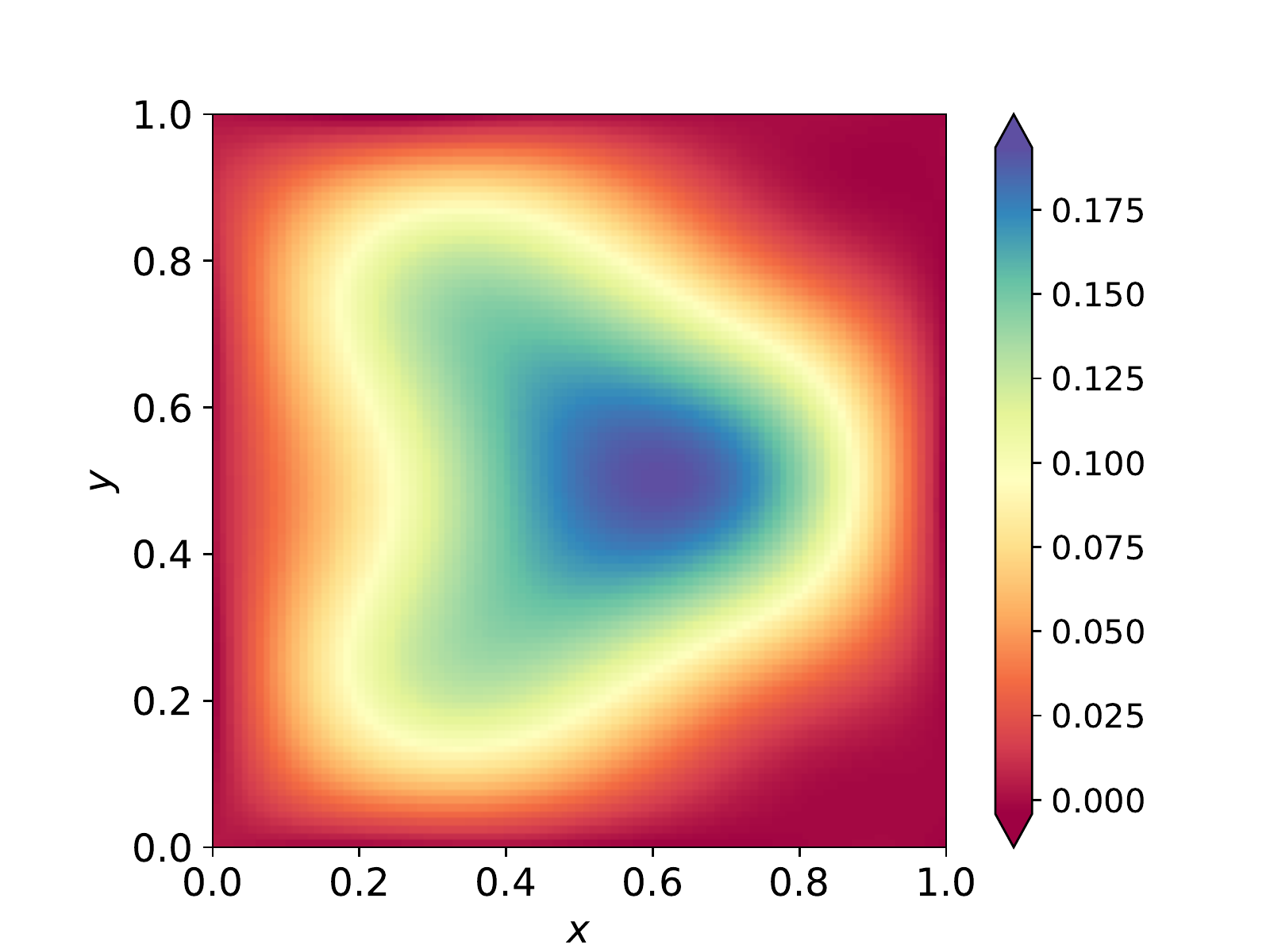}
        \end{minipage}
    }%
    \subfigure[$u-u_{NN}$.]{
    \begin{minipage}[t]{0.3\linewidth}
        \centering
        \includegraphics[width=2in]{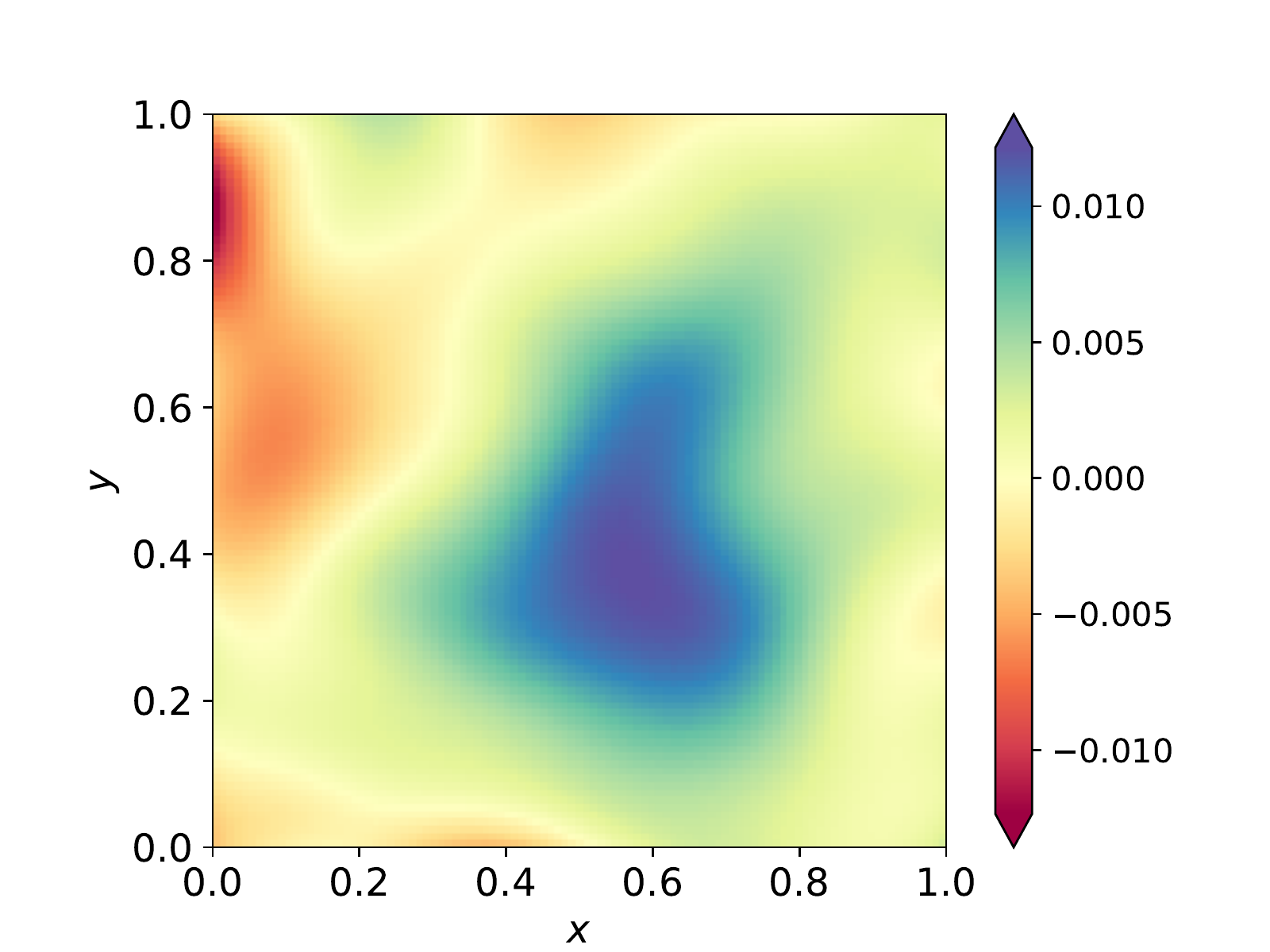}
    \end{minipage}
    }%

    \subfigure[The reference solution $u$.]{
        \begin{minipage}[t]{0.3\linewidth}
            \centering
            \includegraphics[width=2in]{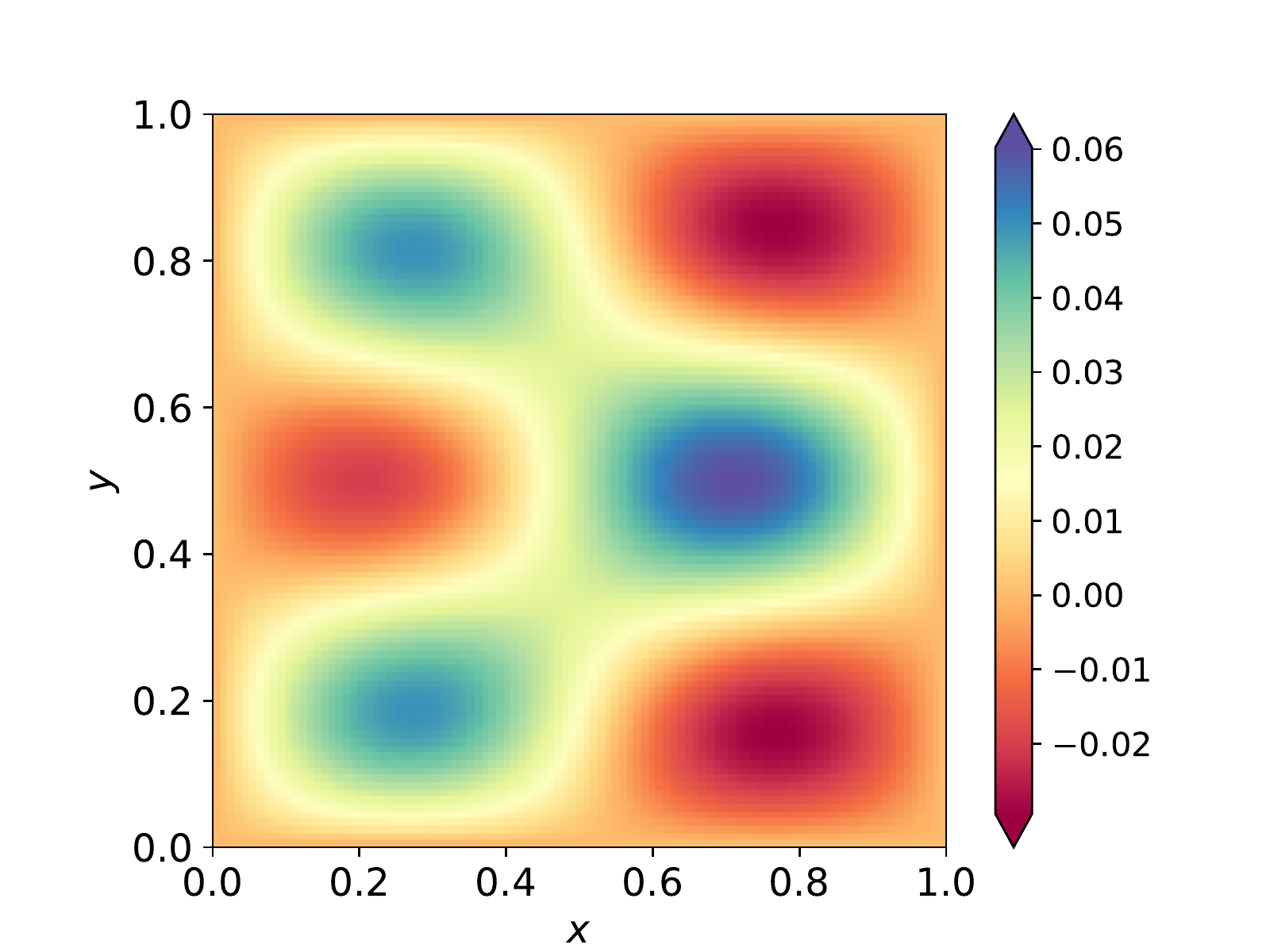}
        \end{minipage}
    }%
    \subfigure[The approximate solution $u_{NN}$.]{
        \begin{minipage}[t]{0.3\linewidth}
            \centering
            \includegraphics[width=2in]{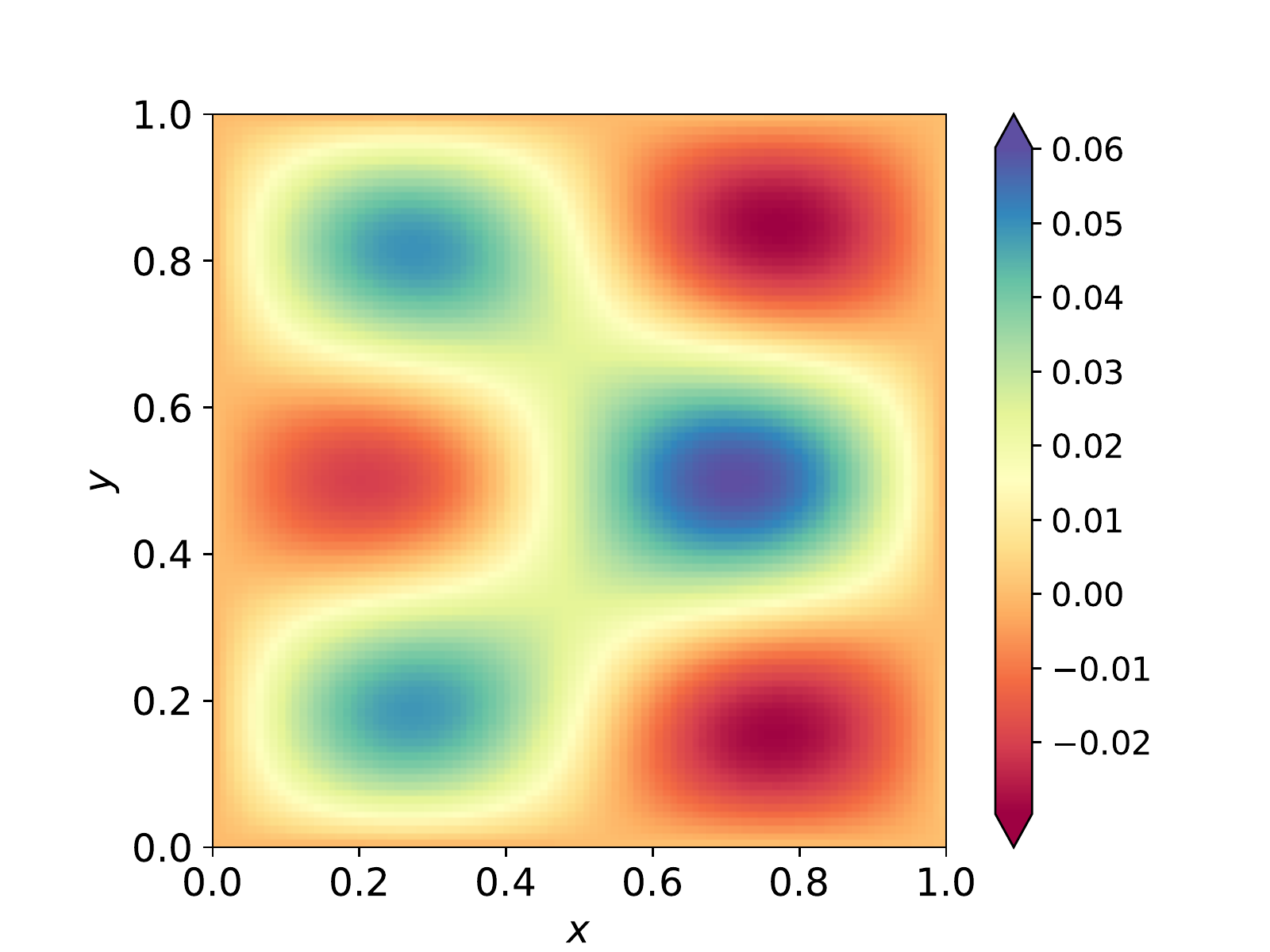}
        \end{minipage}
    }%
    \subfigure[$u-u_{NN}$.]{
    \begin{minipage}[t]{0.3\linewidth}
        \centering
        \includegraphics[width=2in]{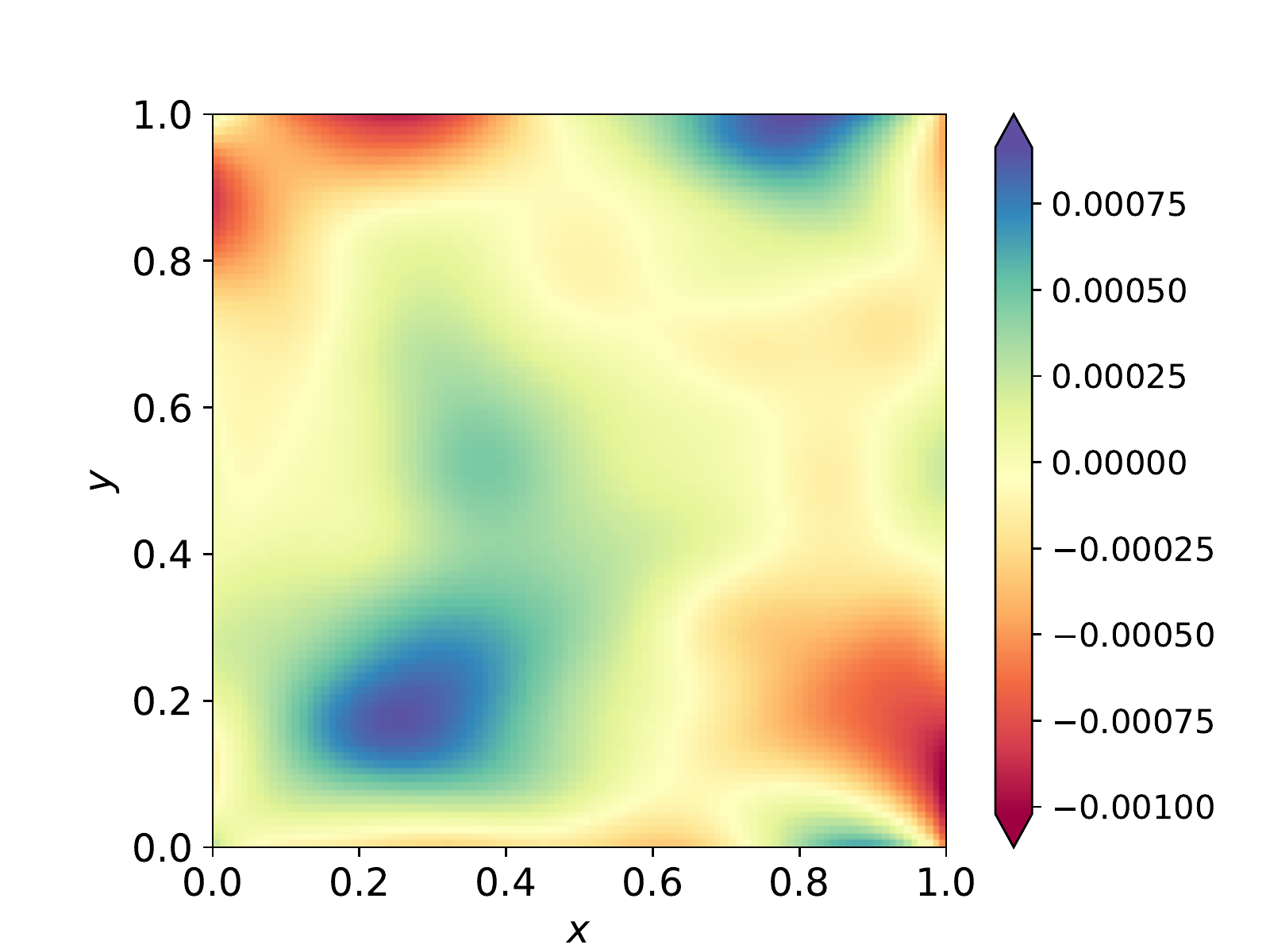}
    \end{minipage}
    }%
    \centering
    \caption{Comparison between the reference solution $u(x,y,t)$ and the approximate solution $u_{NN}(x,y,t;\theta^{*})$ of 2D forward problem with $T=1$. The top panel shows the reference solution $u(x,y,0.02)$, the approximate solution $u_{NN}(x,y,0.02;\theta^{*})$, and the errors $u(x,y,0.02)-u_{NN}(x,y,0.02;\theta^{*})$, respectively. The bottom panel shows the reference solution $u(x,y,1)$, the approximate solution $u_{NN}(x,y,1;\theta^{*})$, and the errors $u(x,y,1)-u_{NN}(x,y,1;\theta^{*})$, respectively.}\label{fig_1}
\end{figure}

In the following, we provide a visual demonstration of the accuracy of the Laplace-fPINNs method in solving the 2D forward problem of equation (\ref{1.1}). Based on the previous numerical results provided in Table~\ref{table_1} and Table~\ref{table_2}, we select a 7-layer fully-connected neural network with 512 neurons per hidden layer for the approximate solution and set $M=4$. Figure~\ref{fig_1} displays the exact solution $u$, the approximate solution $u_{NN}$, and their corresponding residual $u-u_{NN}$. From the residual as shown in Figure~\ref{fig_1} (c), (f), we see that the approximate solution obtained by the Laplace-fPINNs method is in excellent agreement with the ground truth.
\begin{table}[H]
    \centering
    \begin{tabular}{ccccc}
        \toprule
        M                    & 2                 & 4                 & 6                 \\
        \midrule
        relative $l_2$ error & 0.0845$\pm$0.0158 & 0.0348$\pm$0.0103 & 0.5212$\pm$0.1770 \\
        \midrule
    \end{tabular}%
    \caption{Mean and standard deviation of relative $l_2$ error of the
        approximate solution $u_{NN}(x,t;\theta^{*})$ in 2D forward problem with different $M$. The depth and width of the
        neural network is 5 and 256, respectively.}\label{table_1}
\end{table}%

\begin{table}[H]
    \centering
    \begin{tabular}{cccc}
        \toprule
        \multirow{2}{*}{Depth} & \multicolumn{3}{c}{Width}                                         \\
        \cmidrule(lr){2-4}
                               & 128                       & 256               & 512               \\
        \midrule
        5                      & 0.0822$\pm$0.0261         & 0.0348$\pm$0.0103 & 0.0349$\pm$0.0084 \\
        6                      & 0.0994$\pm$0.0388         & 0.0358$\pm$0.0051 & 0.0261$\pm$0.0053 \\
        7                      & 0.1055$\pm$0.0186         & 0.0375$\pm$0.0129 & 0.0204$\pm$0.0062 \\
        \bottomrule
    \end{tabular}
    \caption{Mean and standard derivation of relative $l_2$ error between the approximate
        solution $u_{NN}(x,t;\theta^{*})$ and the reference solution $u$ for different neural network architectures. Here, we fix the $M=4$.}\label{table_2}
\end{table}

\subsubsection{2D forward problem with $T=10$}

\begin{figure}[H]
    \centering
    \subfigure[The reference solution $u$.]{
        \begin{minipage}[t]{0.3\linewidth}
            \centering
            \includegraphics[width=2in]{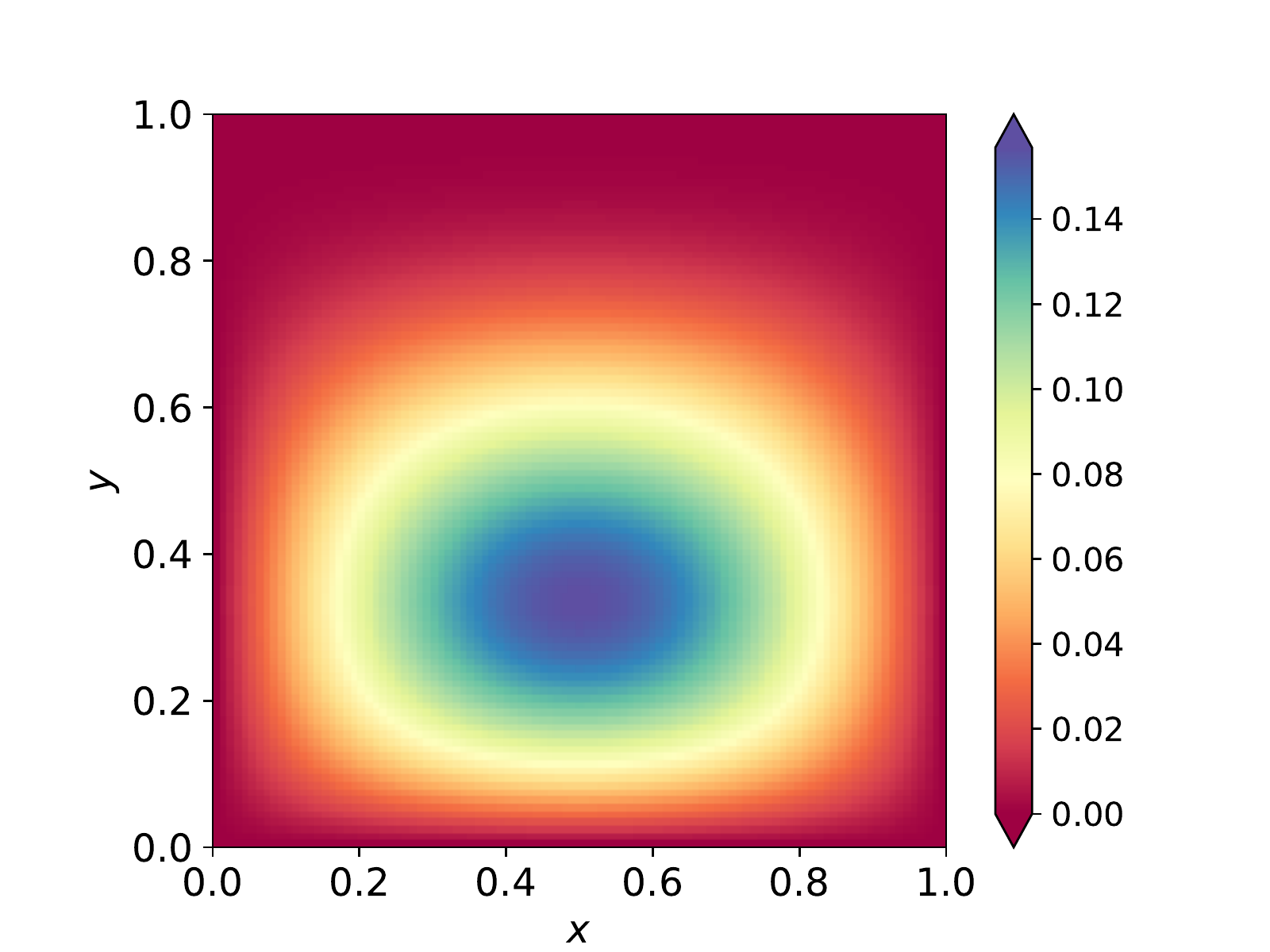}
        \end{minipage}
    }%
    \subfigure[The approximate solution $u_{NN}$.]{
        \begin{minipage}[t]{0.3\linewidth}
            \centering
            \includegraphics[width=2in]{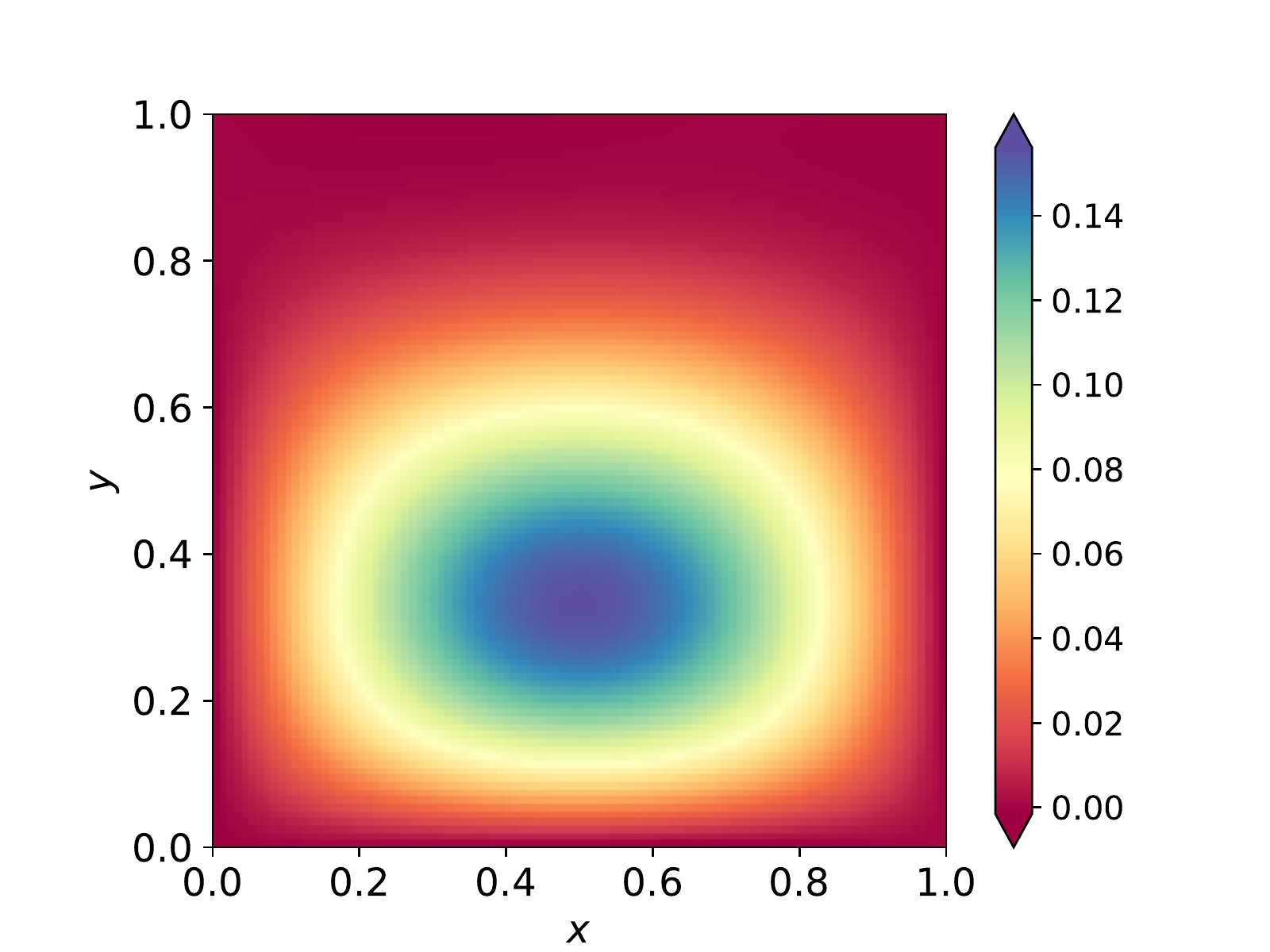}
        \end{minipage}
    }%
    \subfigure[$u-u_{NN}$.]{
    \begin{minipage}[t]{0.3\linewidth}
        \centering
        \includegraphics[width=2in]{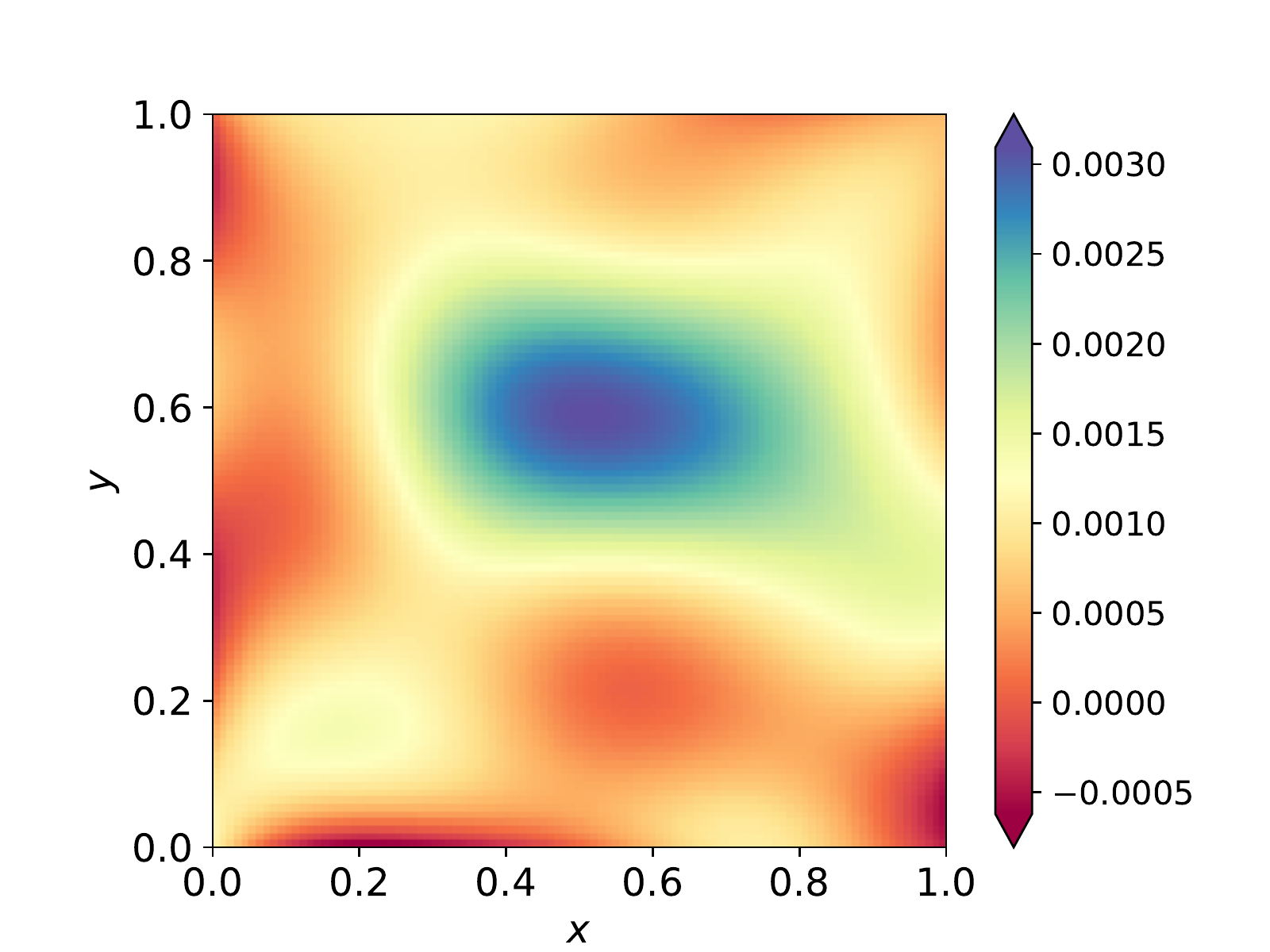}
    \end{minipage}
    }%

    \subfigure[The reference solution $u$.]{
        \begin{minipage}[t]{0.3\linewidth}
            \centering
            \includegraphics[width=2in]{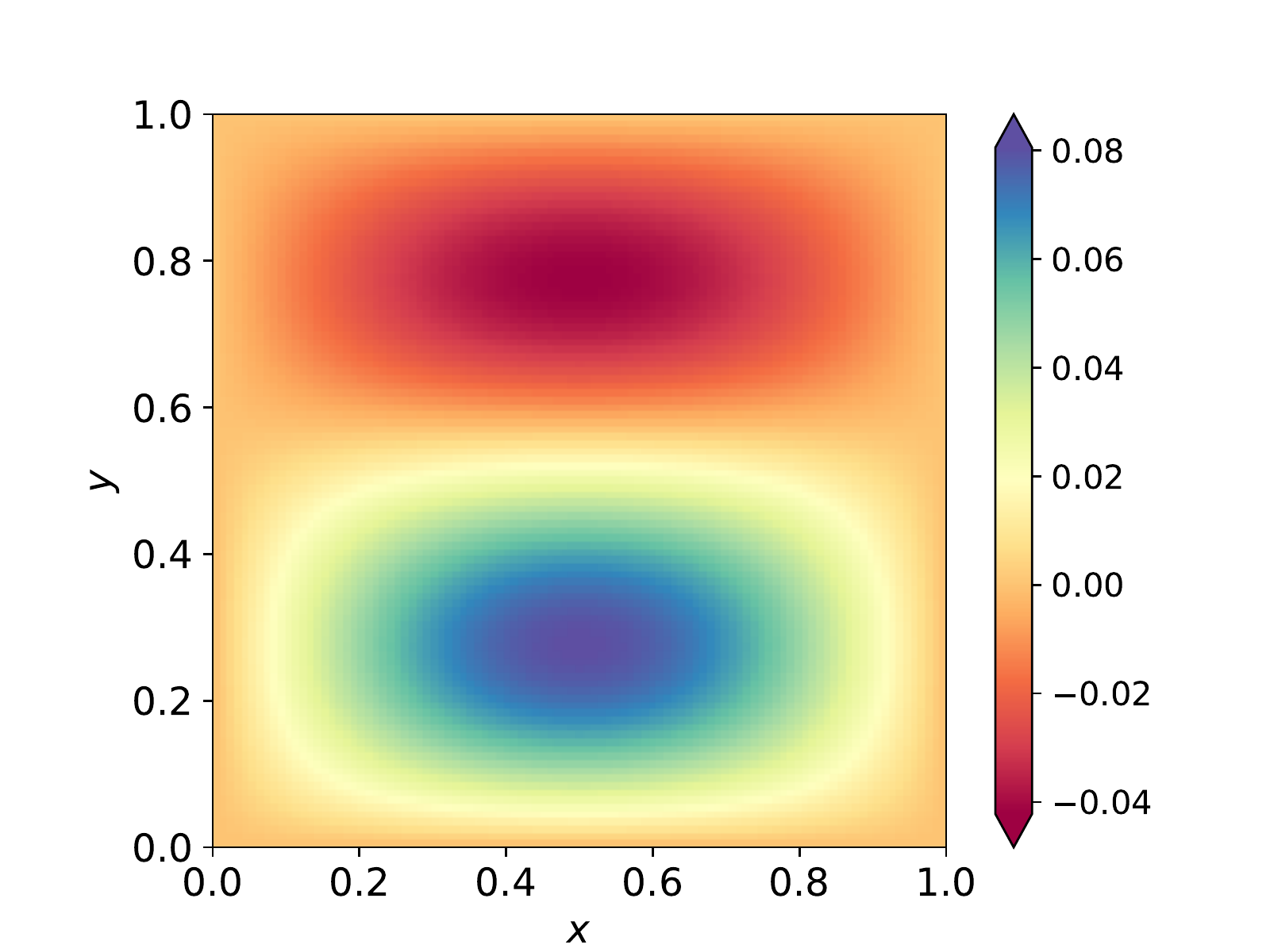}
        \end{minipage}
    }%
    \subfigure[The approximate solution $u_{NN}$.]{
        \begin{minipage}[t]{0.3\linewidth}
            \centering
            \includegraphics[width=2in]{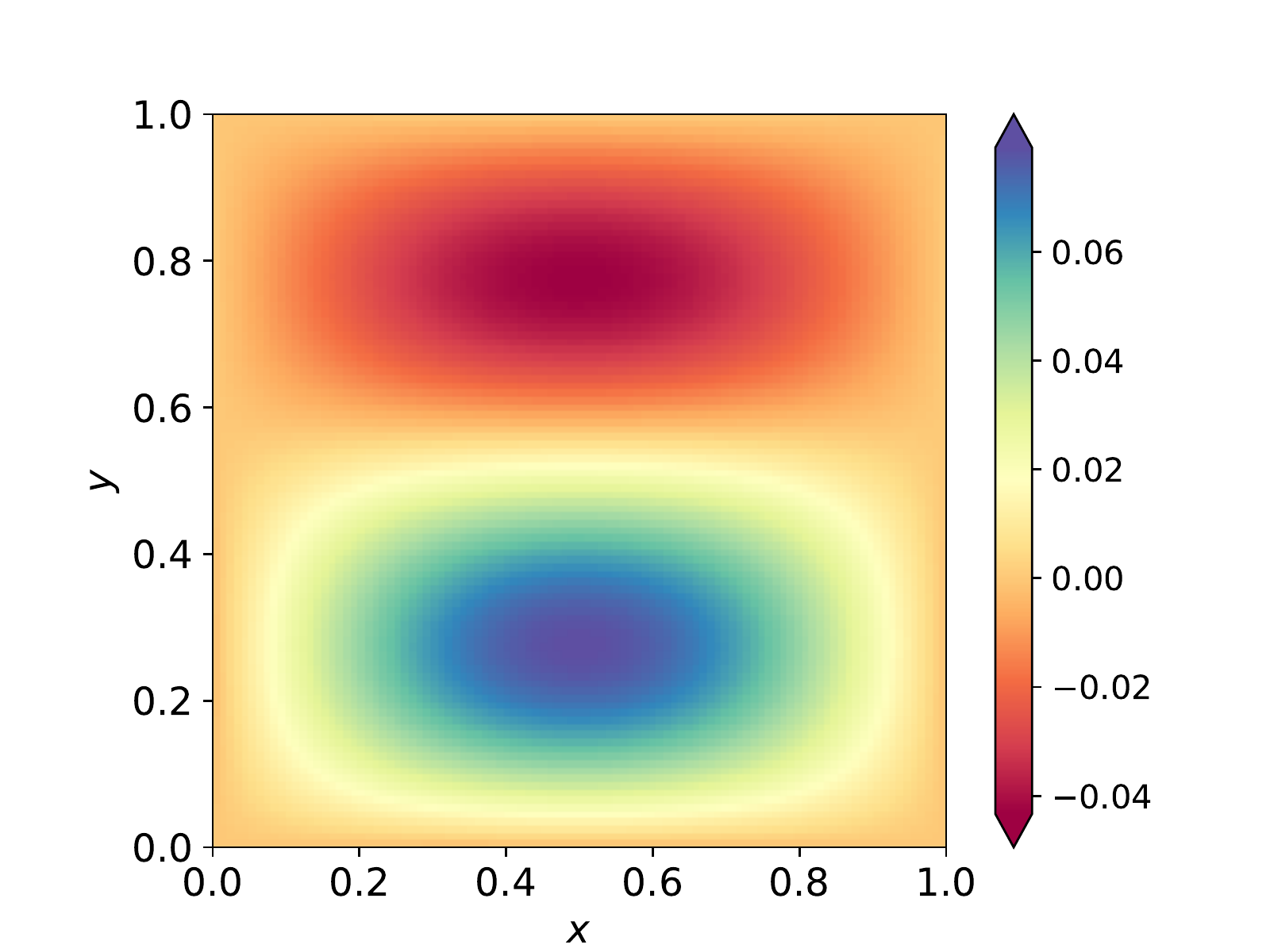}
        \end{minipage}
    }%
    \subfigure[$u-u_{NN}$.]{
    \begin{minipage}[t]{0.3\linewidth}
        \centering
        \includegraphics[width=2in]{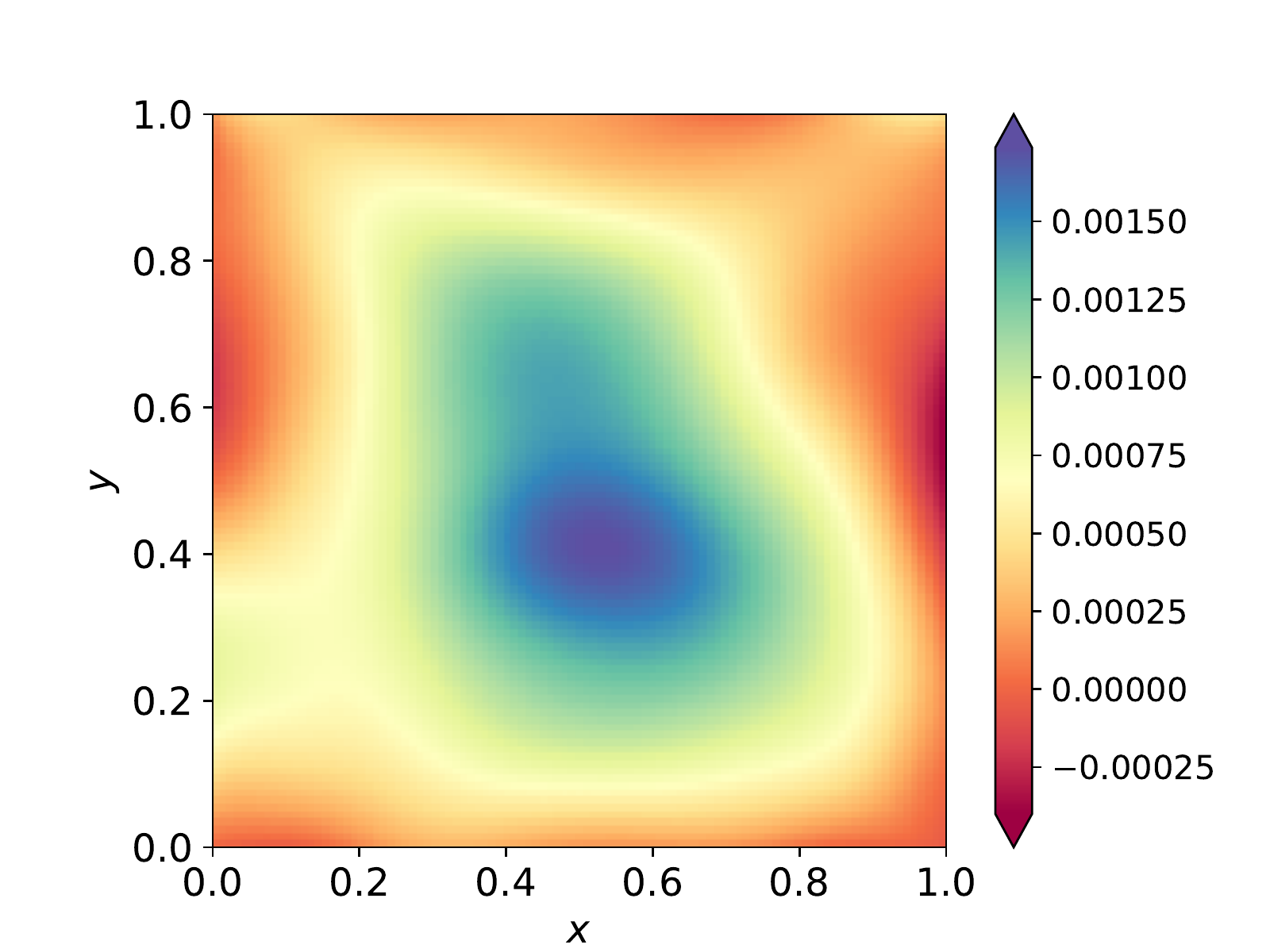}
    \end{minipage}
    }%
    \centering
    \caption{Comparison between the reference solution $u(x,y,t)$ and the approximate solution
    $u_{NN}(x,y,t;\theta^{*})$ of 2D forward problem with $T=10$. The top panel shows the reference solution $u(x,y,0.05)$, the approximate solution $u_{NN}(x,y,0.05;\theta^{*})$, and the errors $u(x,y,0.05)-u_{NN}(x,y,0.05;\theta^{*})$, respectively. The bottom panel shows the reference solution $u(x,y,10)$, the approximate solution $u_{NN}(x,y,10;\theta^{*})$, and the errors $u(x,y,10)-u_{NN}(x,y,10;\theta^{*})$, respectively.}\label{fig_1.1}
\end{figure}

In this part, we consider another example of using the Laplace-fPINNs method
to solve two-dimensional forward problem. Similarly, the diffusion and
reaction coefficients are $a(x,y)=1$, $c(x,y)=0$, respectively. The initial
value and source are given by
\begin{align*}
    u_0(x,y) & =3\sin(\pi x)\sin(\pi y),  \\
    f(x,y)   & =3\sin(\pi x)\sin(2\pi y).
\end{align*}
In this example, we consider the domain $\Omega={[0,1]}^2$ and $T=10$. The neural network architecture, optimizer, and sample point settings are kept the same as those employed in the previous example. As the process of minimizing the loss function for the fPINNs method~\cite{Pang+Lu+Karniadakis-2019}, which approximates the time-fractional derivative via the $L_1$-type finite difference scheme, requires more auxiliary points for large $T$, it can lead to a significant increase in computational cost. However, with the Laplace-fPINNs approach proposed in this work, each iteration's computing cost does not increase with larger $T$. The numerical results, as presented in Figure~\ref{fig_1.1}, illustrate that the Laplace-fPINNs method can accurately solve the forward problem of subdiffusion (\ref{1.1}) even for large $T$.

\subsection{3D forward problem}
\begin{figure}[H]
    \centering
    \subfigure[Exact solution at $t=1,~z=0.5$]{
        \begin{minipage}[t]{0.3\linewidth}
            \centering
            \includegraphics[width=2in]{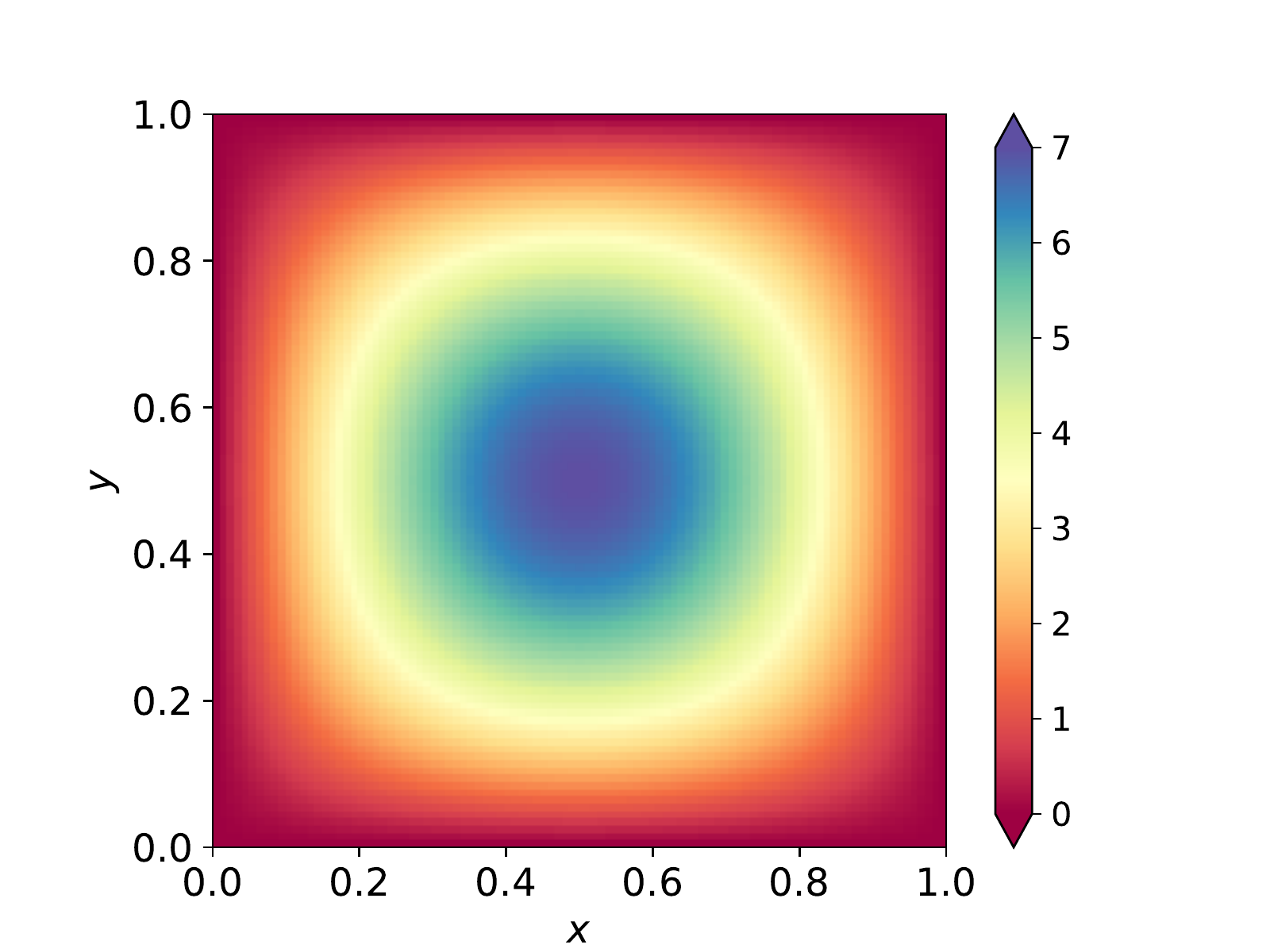}
        \end{minipage}
    }%
    \subfigure[Approximate solution at $t=1,~z=0.5$]{
        \begin{minipage}[t]{0.3\linewidth}
            \centering
            \includegraphics[width=2in]{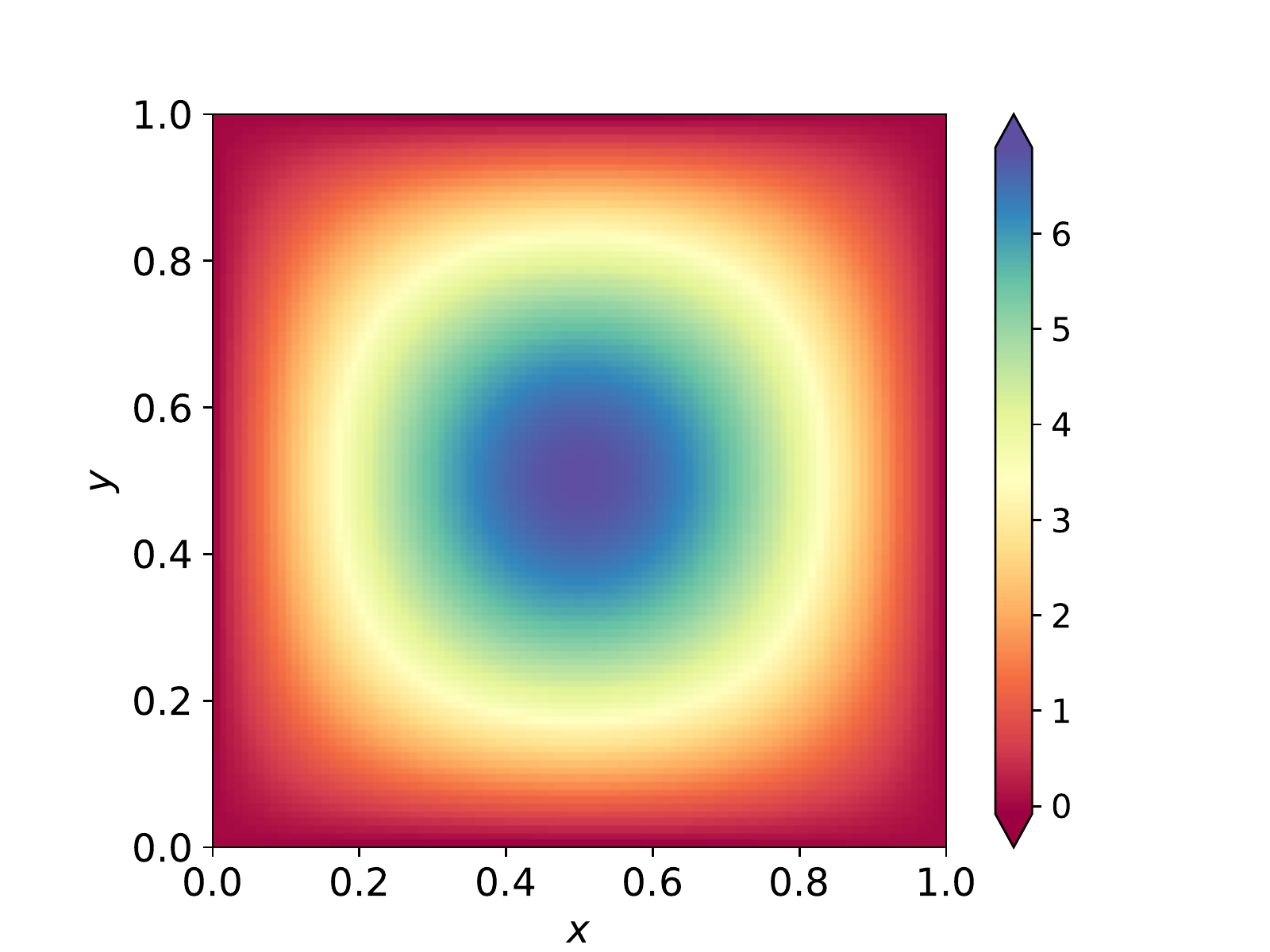}
        \end{minipage}
    }%
    \subfigure[Errors at $t=1,~z=0.5$]{
        \begin{minipage}[t]{0.3\linewidth}
            \centering
            \includegraphics[width=2in]{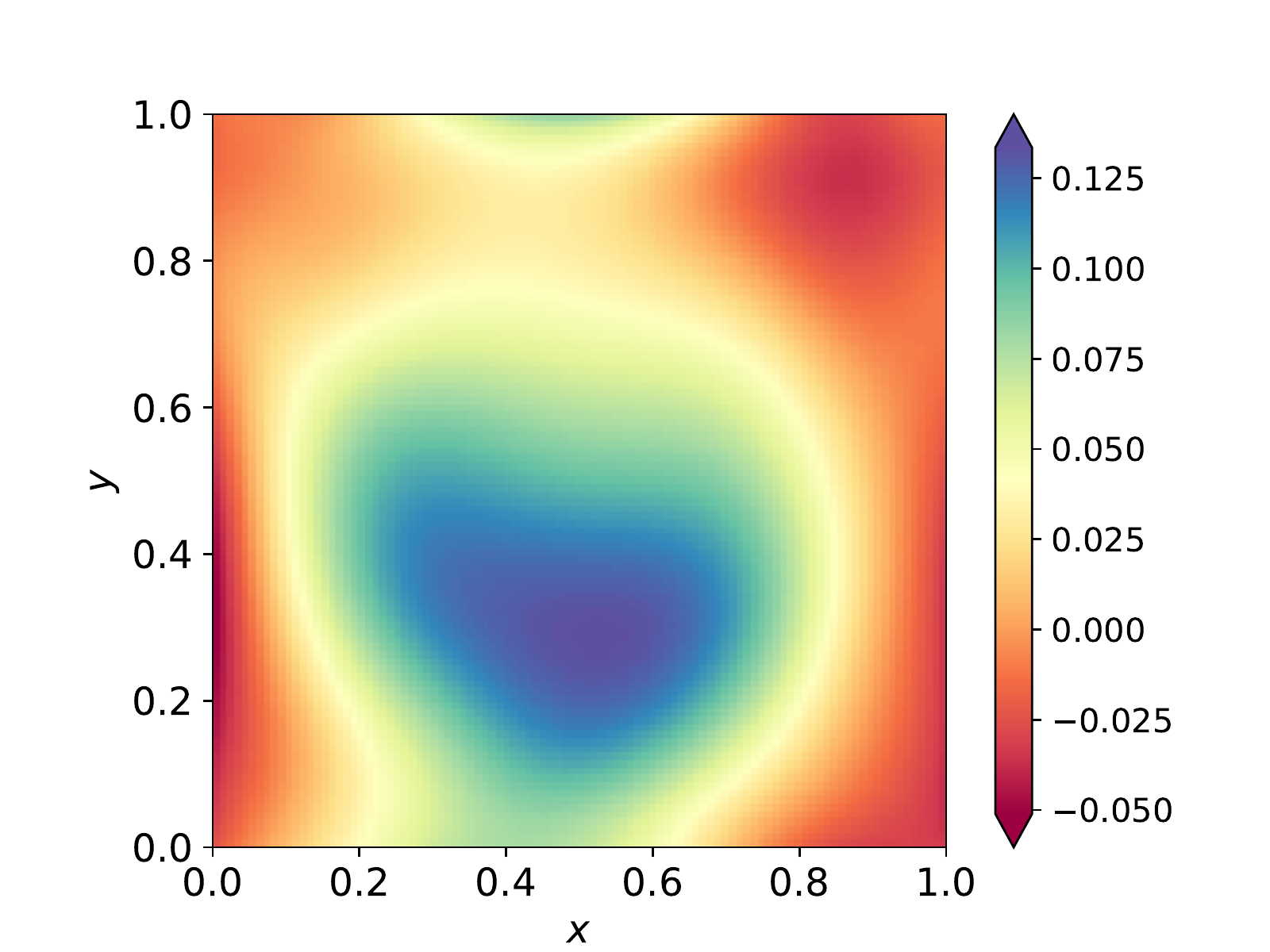}
        \end{minipage}
    }%

    \subfigure[Exact solution at $t=0.5,~y=0.5$]{
        \begin{minipage}[t]{0.3\linewidth}
            \centering
            \includegraphics[width=2in]{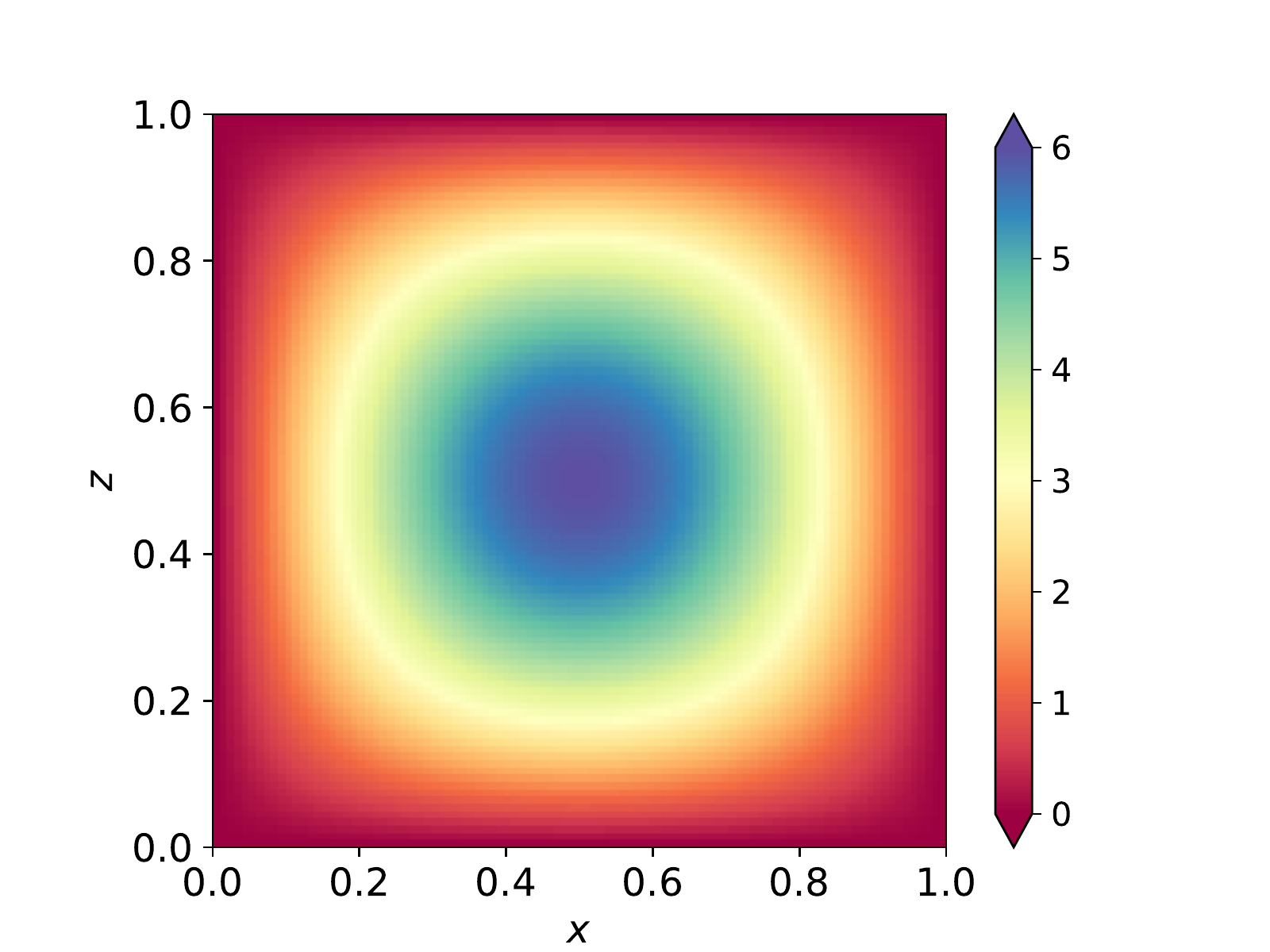}
        \end{minipage}
    }%
    \subfigure[Approximate solution at $t=0.5,~y=0.5$]{
        \begin{minipage}[t]{0.3\linewidth}
            \centering
            \includegraphics[width=2in]{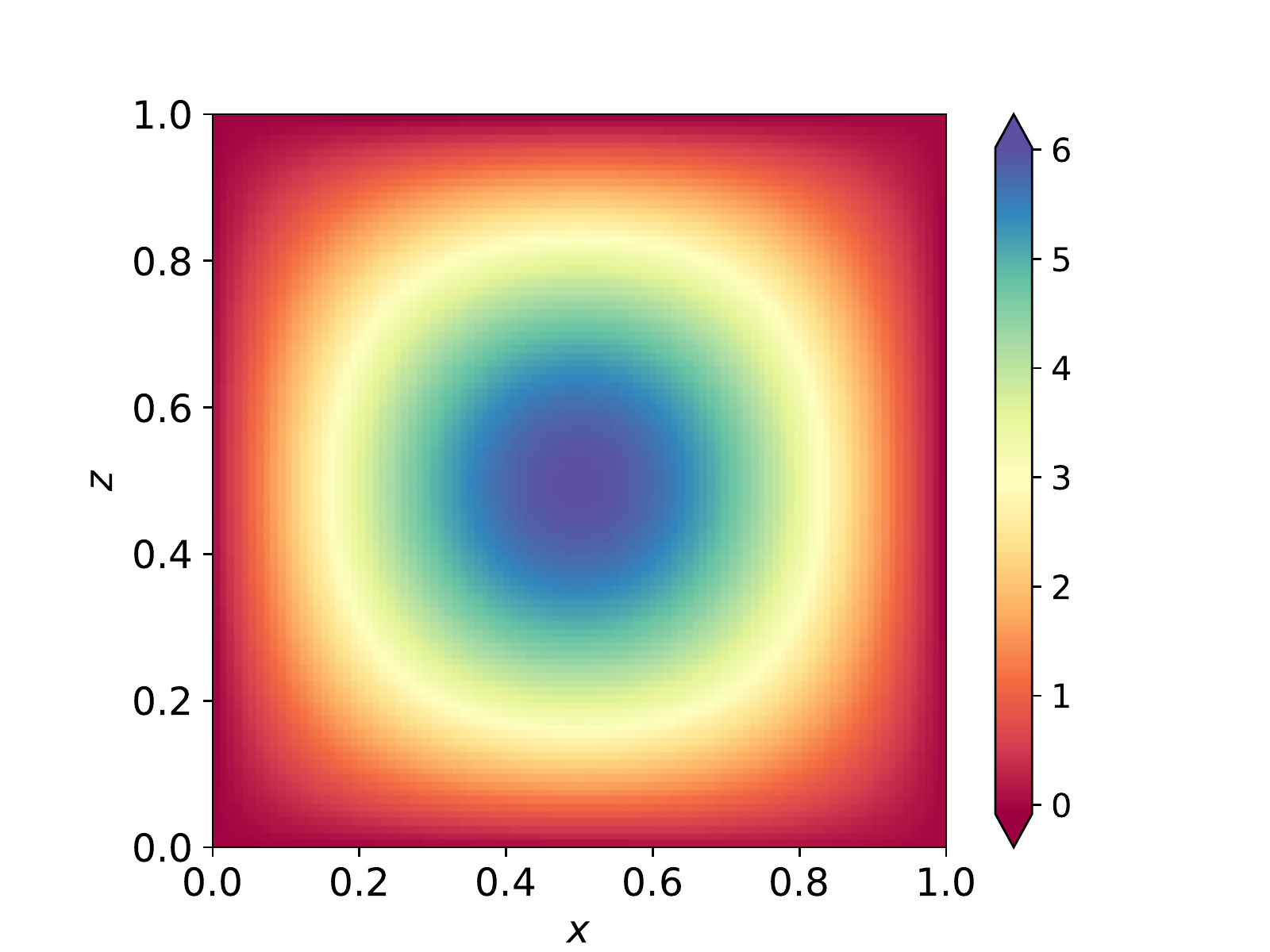}
        \end{minipage}
    }%
    \subfigure[Errors at $t=0.5,~y=0.5$]{
        \begin{minipage}[t]{0.3\linewidth}
            \centering
            \includegraphics[width=2in]{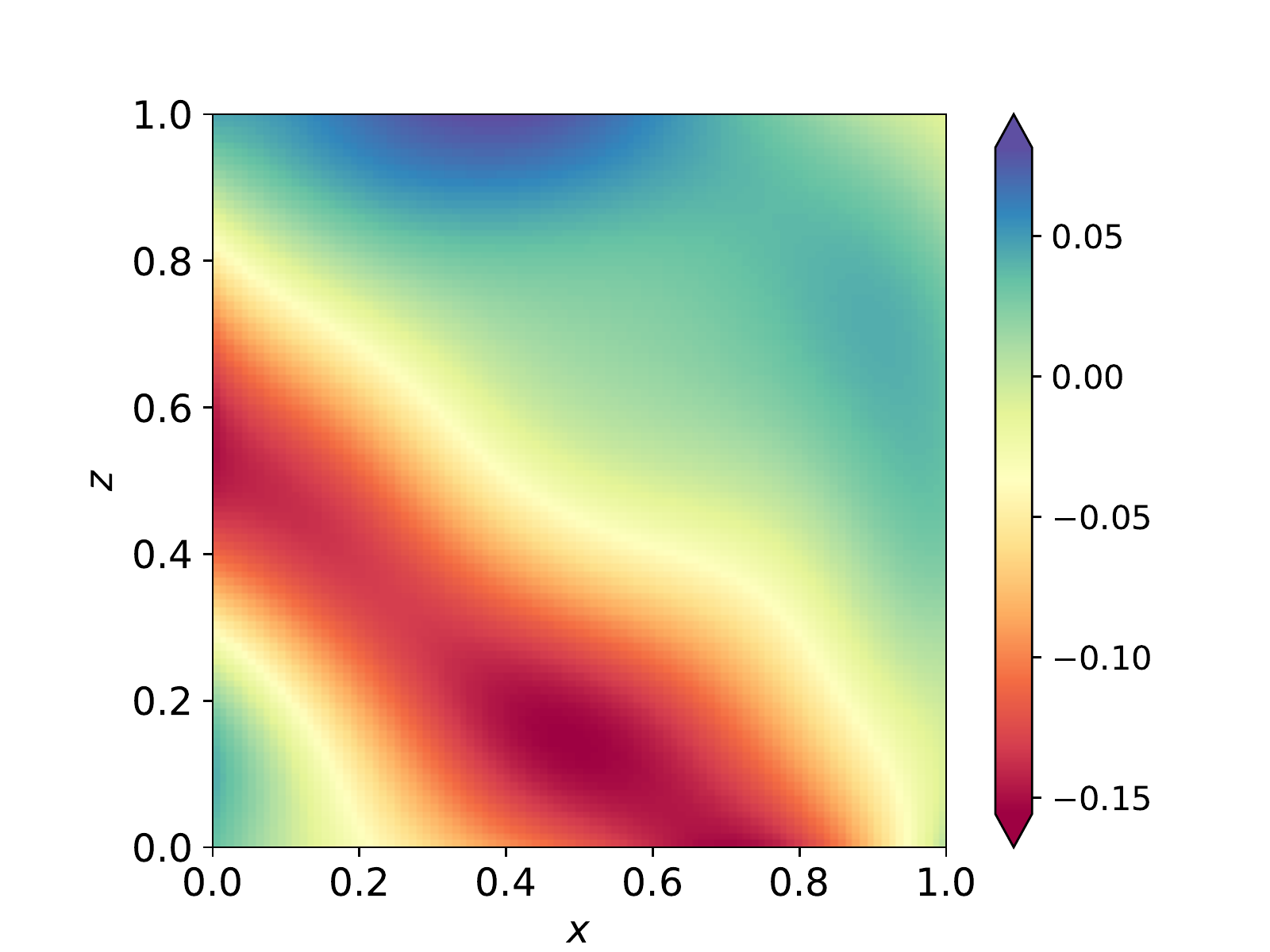}
        \end{minipage}
    }%
    \centering
    \caption{Comparison between the reference solution $u$ and the approximate solution
    $u_{NN}(x,y,z,t;\theta^{*})$ of 3D forward problem. The top panel shows the reference solution $u(x,y,0.5,1)$, the approximate solution $u_{NN}(x,y,0.5,1;\theta^{*})$, and the errors $u(x,y,0.5,1)-u_{NN}(x,y,0.5,1;\theta^{*})$, respectively. The bottom panel shows the reference solution $u(x,0.5,z,0.5)$, the approximate solution $u_{NN}(x,0.5,z,0.5;\theta^{*})$, and the errors $u(x,0.5,z,0.5)-u_{NN}(x,0.5,z,0.5;\theta^{*})$, respectively.}\label{fig_2}
\end{figure}

In this part, we study the accuracy of using the Laplace-fPINNs method to
solve a 3D forward problem in equation (\ref{1.1}). We take the domain
$\Omega={[0,1]}^3$, $T=1$, the diffusion coefficient $a(x,y,z)=1$, and the reaction coefficient $c(x,y,z)=0$. A manufactured solution is taken as
$u(x,y,z,t)=(2t+5)\sin(\pi x)\sin(\pi y)\sin(\pi z)$, and the corresponding
initial value and source are given by
\begin{align*}
    u_0(x,y,z) & =5\sin(\pi x)\sin(\pi y)\sin(\pi z),                                                                        \\
    f(x,y,z,t) & =[\frac{2}{\Gamma(2-\alpha)}t^{1-\alpha}+(2t+5)(3\pi^2a(x,y,z)-c(x,y,z))]\sin(\pi x)\sin(\pi y)\sin(\pi z).
\end{align*}
We choose a fully-connected
neural network with 5 layers and 256 neurons per layer and set the iteration number
to $1.5\times 10^{5}$. The numerical results obtained from the trained model are illustrated in Figure~\ref{fig_2}. The top panel of Figure~\ref{fig_2} displays the comparison between the exact solution $u$ and the approximate solution $u_{NN}$, with fixed $t=1,~z=0.5$. Similarly, the bottom panel of Figure~\ref{fig_2} shows the comparision between the exact solution $u$ and the approximate solution $u_{NN}$, where we fix $t=0.5,~y=0.5$. From Figure~\ref{fig_2}, we can find that the Laplace-fPINNs method can achieve the higher accuracy even for the 3D forward problem in equation (\ref{1.1}).

\subsection{Inverse problem}
In this part, we study the performance of using the Laplace-fPINNs to solve an inverse problem of the time-fractional diffusion equation (\ref{1.1}). We consider a fabricated solution $u(x,y,z,t)=(3t+5)\sin(2\pi x)\sin(\pi y)\sin(\pi z)$ and choose domain $\Omega={[0,1]}^3$, $T=1$, respectively. According to the definition of Caputo derivative, the time-fractional derivative can be computed analytically, and thus the initial value and source are given by
\begin{align*}
    u_0(x,y,z) & =5\sin(2\pi x)\sin(\pi y)\sin(\pi z),                                      \\
    f(x,y,z,t) & =[\frac{3}{\Gamma(2-\alpha)}t^{1-\alpha}+(3t+5)(6\pi^2 a(x,y,z)-c(x,y,z))]
    \sin(2\pi x)\sin(\pi y)\sin(\pi z),
\end{align*}
where the diffusion coefficient $a(x,y,z)=0.5+e^{-(x+y+z)}$, the reaction coefficient $c(x,y,z)=0$. In this paper, we consider the inverse problem of identifying the diffusion coefficient $a$ and solution $u$ from additional measurements of $u$, given that the initial value $u_0$, source $f$, fractional order $\alpha$, and reaction coefficient $c=0$ are known. Specifically, the extra measurements are given by
\begin{align*}
    h(x,y,z,t)=u(x,y,z,t)+\epsilon u(x,y,z,t),\quad (x,y,z)\in\omega\subset\Omega,~t\in(0,T),
\end{align*}
where we take $\epsilon=0.001$, $\omega={[0.3,0.7]}^3$.

To cope with the inverse problem considered in this paper, we use neural networks to parameterize both the solution $u$ and the diffusion coefficient $a$. Specifically, we parameterize the solution $\tilde{u}$ in equation (\ref{2.5}) using a 5-layer fully-connected neural network with 256 neurons per layer, denoted as $\tilde{u}_{NN}(x,y,z,s;\theta_1)$. Additionally, we represent the unknown diffusion coefficient $a$ as a fully-connected neural network with 4 layers and 64 neurons per layer, denoted as $a_{NN}(x,y,z;\theta_2)$. In addition, we assume that the boundary of the diffusion coefficient is known. The parameters $(\theta_1,\theta_2)$ of the approximate solution $\tilde{u}_{NN}(x,y,z,s;\theta_1)$ and the approximate
diffusion coefficient $a_{NN}(x,y,z;\theta_2)$ are optimized by
minimizing the following loss function
\begin{align}\label{3.1}
    (\theta_{1}^{*}, \theta_2^{*}) & =\arg\min L^{lp}(\theta_{1}, \theta_2)                                       \\ \nonumber
                                   & =\arg\min{\{w_{eq}L_{eq}^{lp}(\theta_1,\theta_2)+w_{bd}L_{bd}^{lp}(\theta_1)
    +w_{obs}L_{obs}^{lp}(\theta_1)+w_{prior}L_{prior}^{lp}(\theta_2)\}},
\end{align}
where
\begin{align*}
    L_{eq}^{lp}(\theta_1,\theta_2) & =\frac{1}{N_r}\sum_{i=1}^{N_r}|r_{NN}^{ip}(x_{r}^{i},y_{r}^{i},z_{r}^{i},s_{r}^{i};\theta_1, \theta_2)|^2,                                                                \\
    L_{obs}^{lp}(\theta_1)         & =\frac{1}{N_{obs}}\sum_{i=1}^{N_{obs}}|\tilde{u}_{NN}(x_{obs}^{i}, y_{obs}^{i}, z_{obs}^{i}, s_{obs}^{i};\theta_1)-h(x_{obs}^{i},y_{obs}^{i},z_{obs}^{i},s_{obs}^{i})|^2, \\
    L_{prior}^{lp}(\theta_2)       & =\frac{1}{N_{prior}}\sum_{i=1}^{N_{prior}}|a_{NN}(x_{prior}^{i},y_{prior}^{i},z_{prior}^{i})-a(x_{prior}^{i},y_{prior}^{i},z_{prior}^{i})|^2.
\end{align*}
The $L_{bd}^{lp}(\theta_1)$ defined by (\ref{2.6}) and the $r_{NN}^{ip}(x_r,y_r,z_r,s_r;\theta_1,\theta_2)$ given by
\begin{align*}
    r_{NN}^{ip}(x,y,z,s;\theta_1,\theta_2) & :=s^{\alpha}\tilde{u}_{NN}(x,y,z,s;\theta_1)-\nabla\cdot(a_{NN}(x,y,z;\theta_2)\nabla \tilde{u}_{NN}(x,y,z,s;\theta_1)) \\
                                           & -c(x,y,z)\tilde{u}_{NN}(x,y,z,s;\theta_1)-s^{\alpha-1}u_0(x,y,z)-\tilde{f}(x,y,z,s).
\end{align*}
The $N_{obs}$ and $N_{prior}$ represent the batch sizes of the observation points corresponding to the measurement data, and the known boundary points for the diffusion coefficient $a(x,y,z)$, respectively. For each iteration, we randomly generate 1000 residual points and $400^3$ boundary points from the uniform distribution on $\Omega$ and $\partial\Omega$, respectively. Furthermore, during training, we fix the observation points $N_{obs}=31^3$ and the prior points $N_{prior}=51^3$. The weight coefficients in (\ref{3.1}) are taken as $w_{eq}=1$, $w_{bd}=2000$, $w_{obs}=1000$, $w_{prior}=100$.

\begin{figure}[H]
    \centering
    \subfigure[Exact $a(x,y,z)$ at $z=0.8$]{
        \begin{minipage}[t]{0.3 \linewidth}
            \centering
            \includegraphics[width=2in]{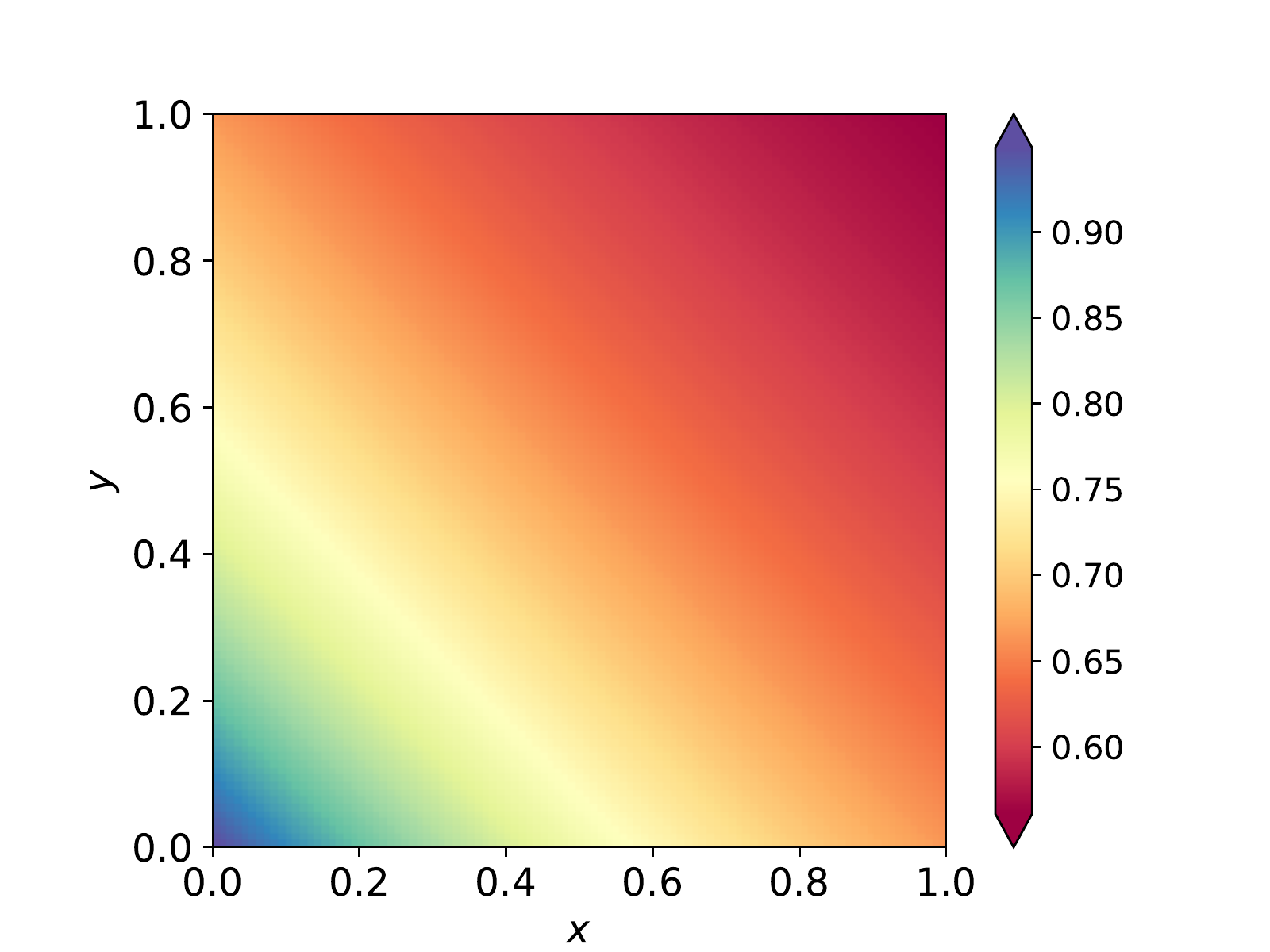}
        \end{minipage}
    }%
    \subfigure[Approximate $a_{NN}(x,y,z;\theta_{2}^{*})$ at $z=0.8$]{
        \begin{minipage}[t]{0.3\linewidth}
            \centering
            \includegraphics[width=2in]{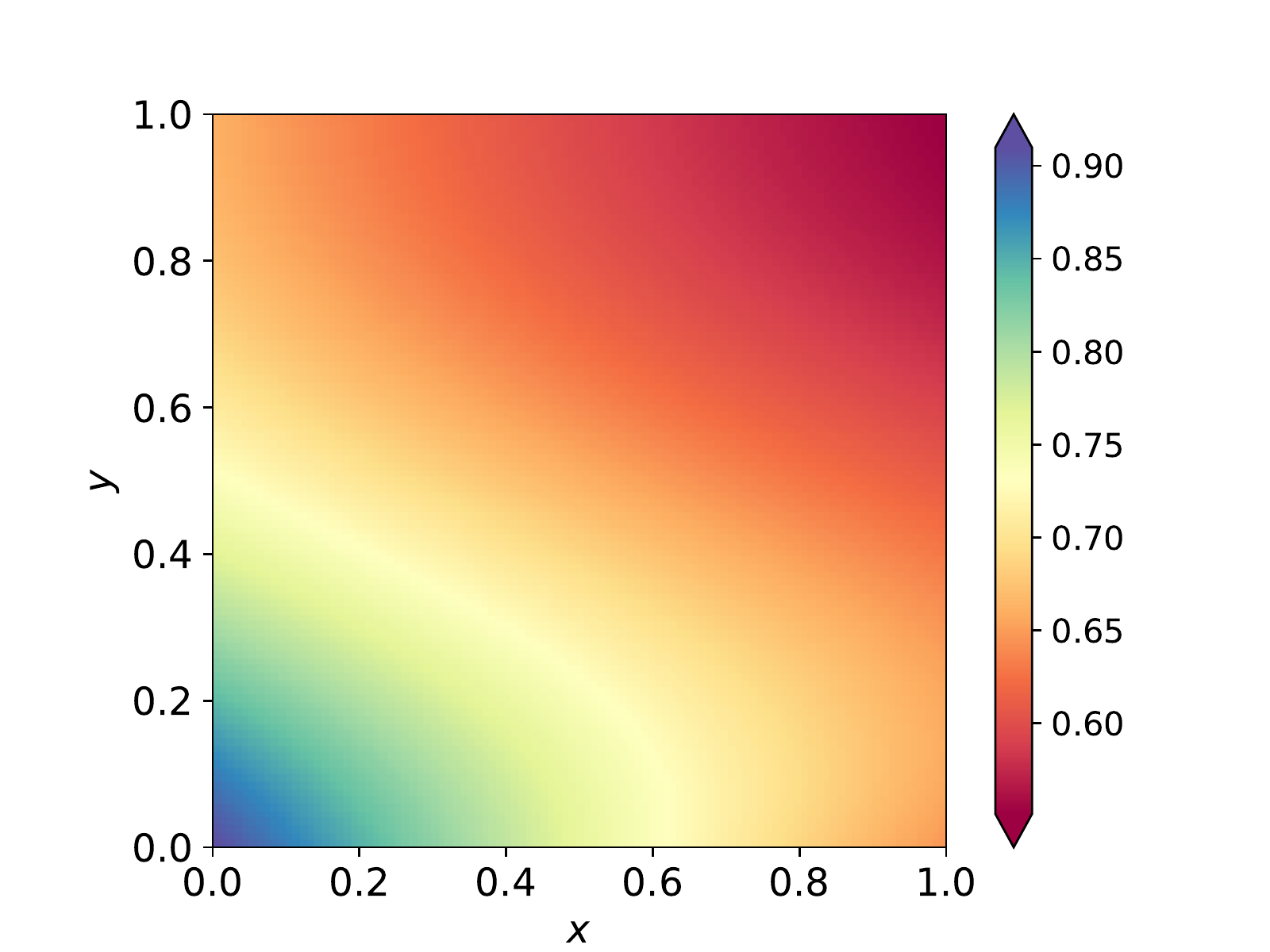}
        \end{minipage}
    }%
    \subfigure[$a(x,y,0.8)-a_{NN}(x,y,0.8;\theta_{2}^{*})$]{
    \begin{minipage}[t]{0.3\linewidth}
        \centering
        \includegraphics[width=2in]{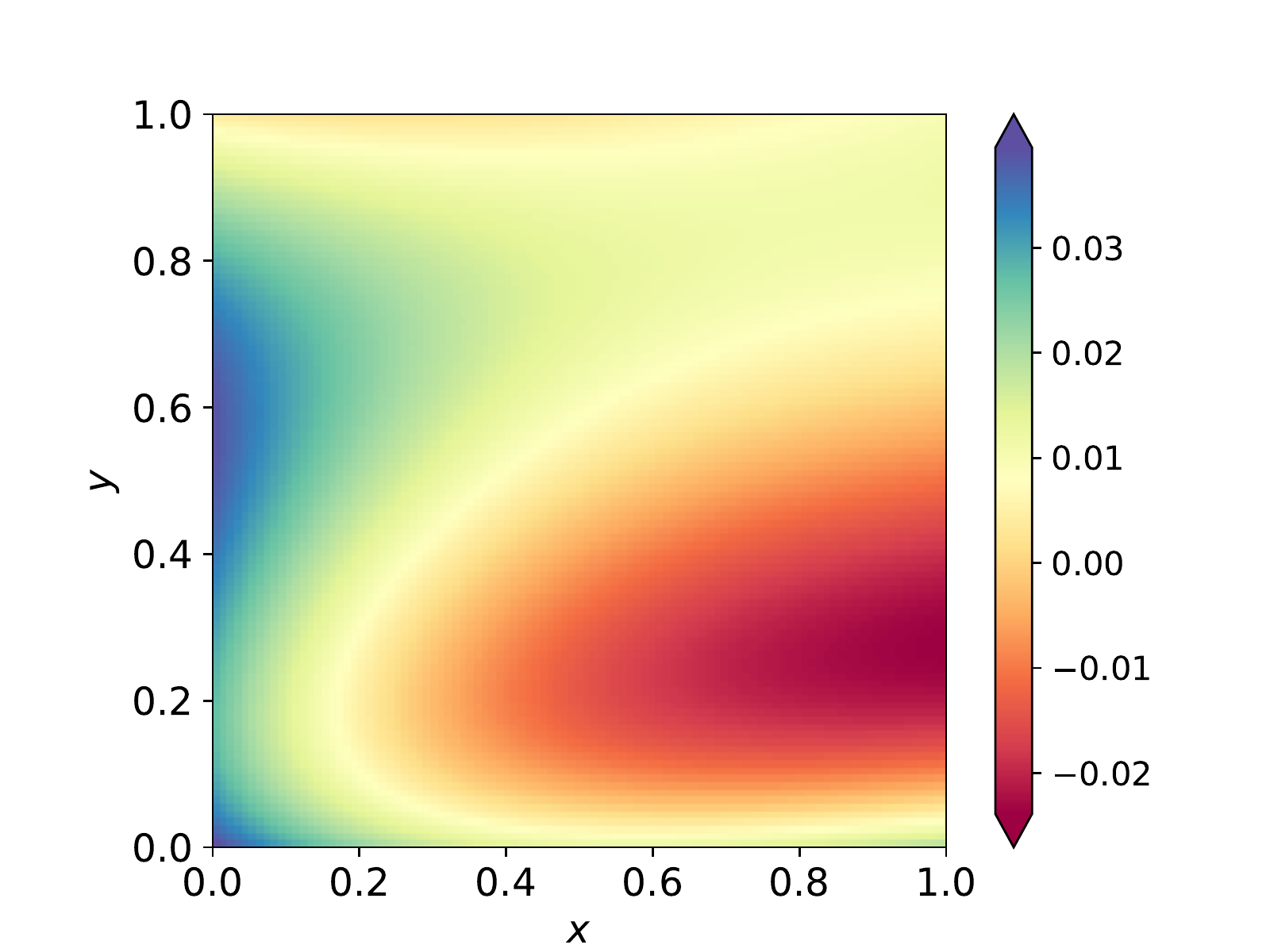}
    \end{minipage}
    }%

    \subfigure[Exact $u(x,y,z,t)$ at $z=0.8$, $t=1$]{
        \begin{minipage}[t]{0.3\linewidth}
            \centering
            \includegraphics[width=2in]{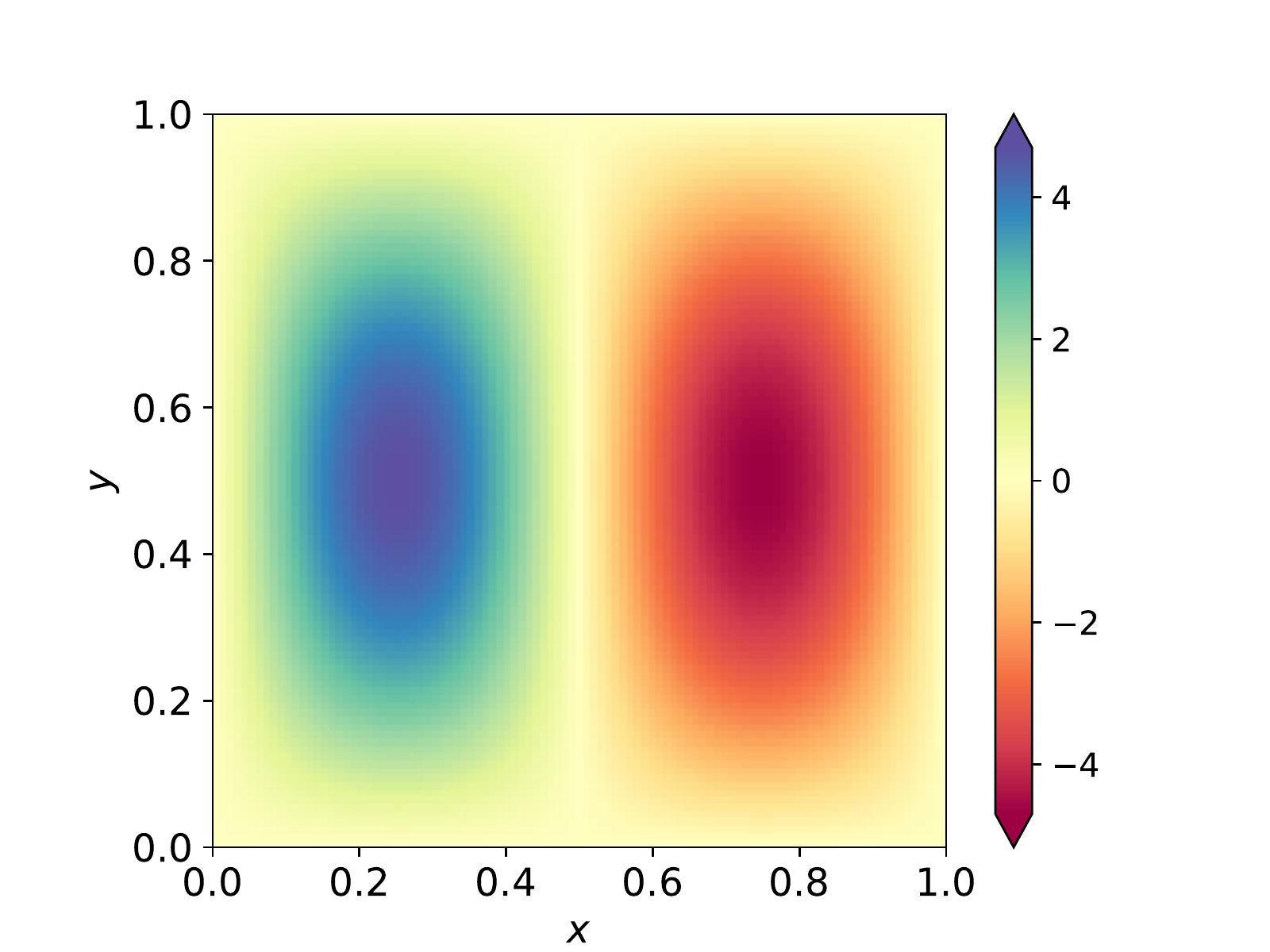}
        \end{minipage}
    }%
    \subfigure[$u_{NN}(x,y,z,t;\theta_{1}^{*})$ at $z=0.8$, $t=1$]{
        \begin{minipage}[t]{0.3\linewidth}
            \centering
            \includegraphics[width=2in]{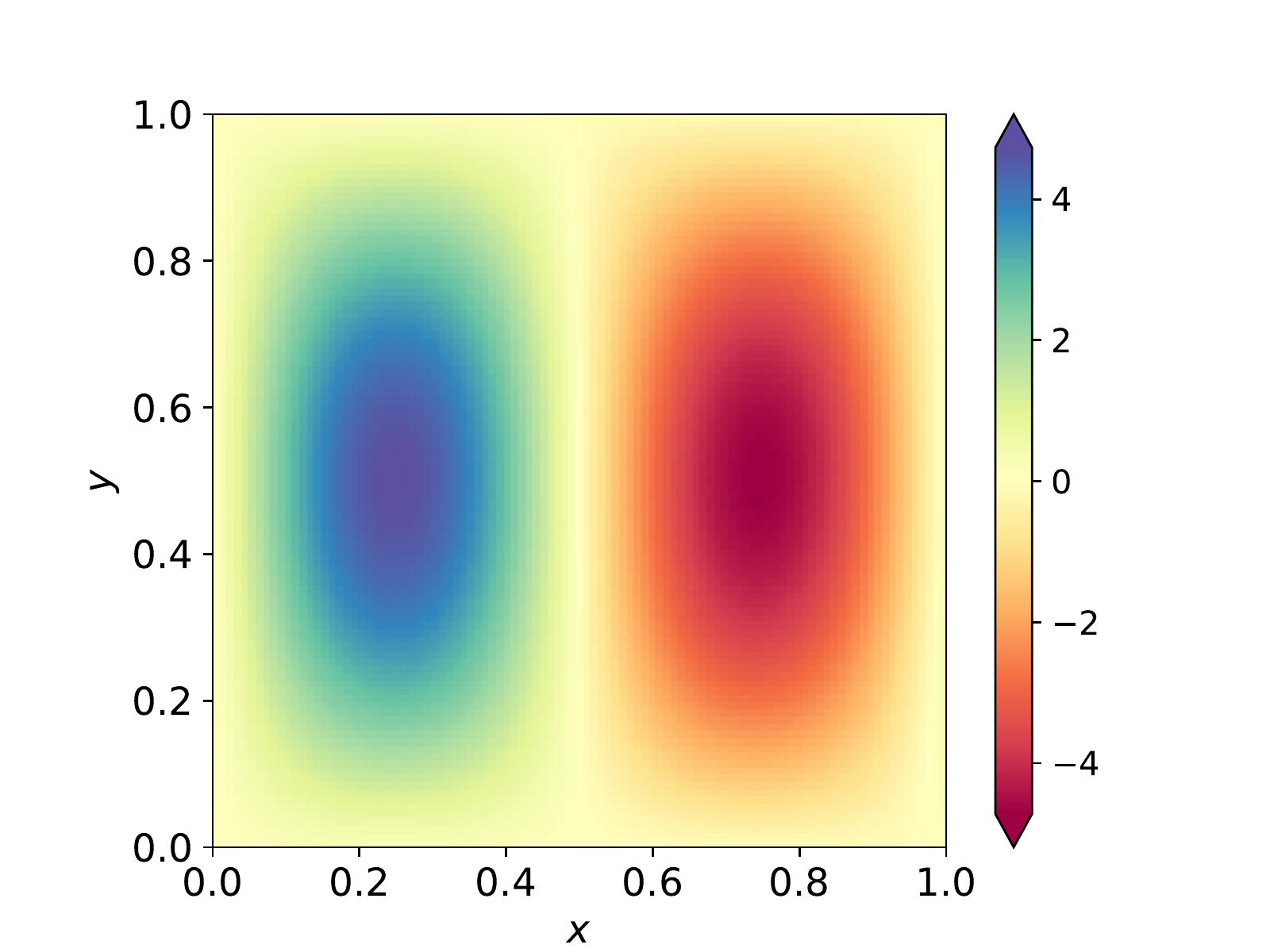}
        \end{minipage}
    }%
    \subfigure[$u(x,y,0.8,1)-u_{NN}(x,y,0.8,1;\theta_{1}^{*})$]{
    \begin{minipage}[t]{0.3\linewidth}
        \centering
        \includegraphics[width=2in]{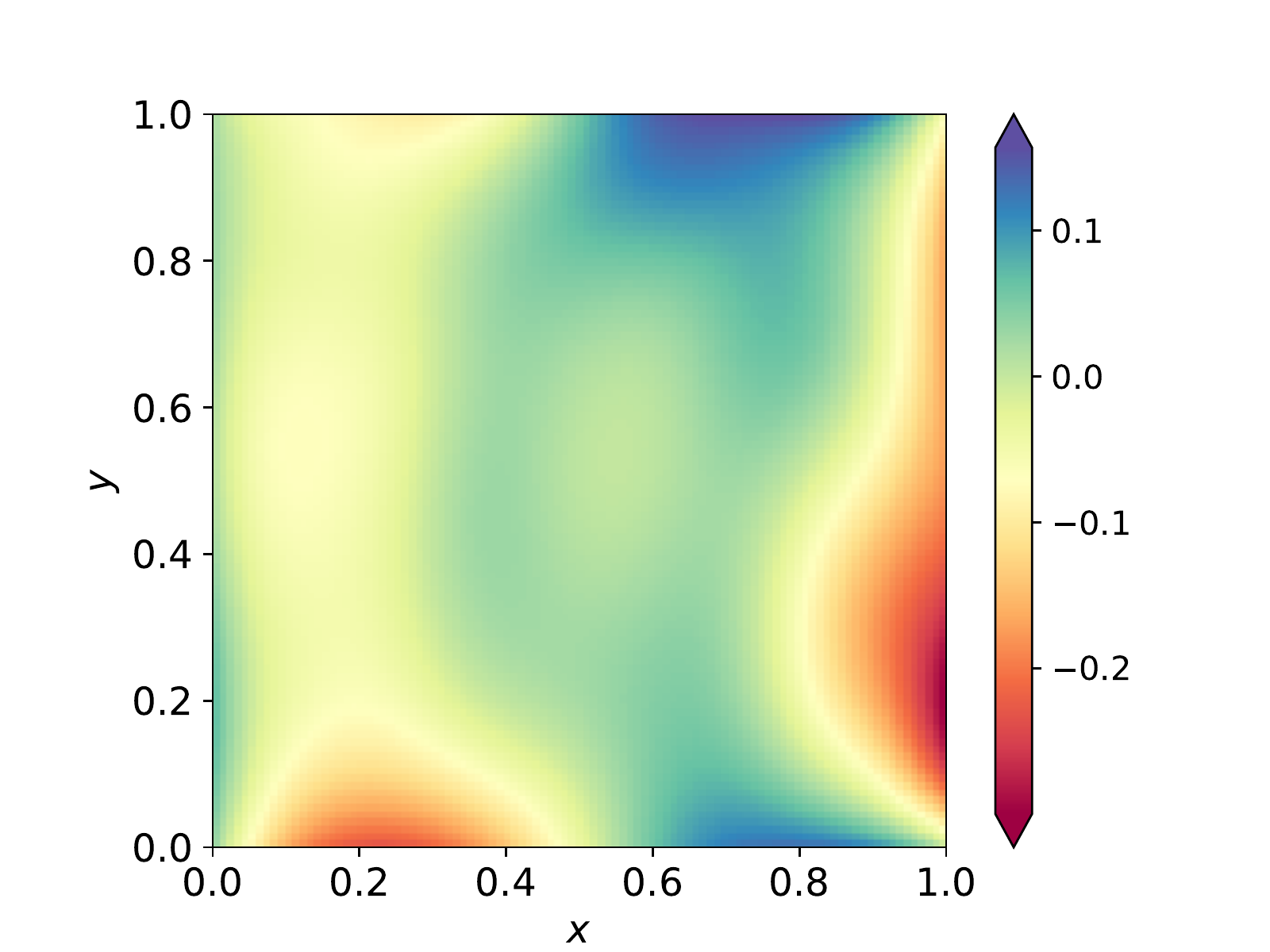}
    \end{minipage}
    }%
    \centering
    \caption{The top panel shows the comparison between the exact diffusion
        coefficient $a(x,y,z)$ and the approximate diffusion coefficient $a_{NN}(x,y,z;\theta_{2}^{*})$ at $z=0.8$. The bottom panel shows comparison between the reference solution $u(x,y,z,t)$ and the approximate solution $u_{NN}(x,y,z,t;\theta_1^{*})$ at $z=0.8$, $t=1$.}\label{fig_3}
\end{figure}

After $1.5\times 10^{5}$ epochs of training, the numerical results obtained
from the trained model are shown in Figure~\ref{fig_3}. Specifically, Figure~\ref{fig_3} (b) displays the reconstructed diffusion coefficient at $z=0.8$, which agrees with the exact diffusion coefficient shown in Figure~\ref{fig_3} (a). Additionally, Figure~\ref{fig_3} (c) presents the difference between the reconstructed diffusion coefficient $a_{NN}(x,y,z;\theta_{2}^{*})$ and the exact diffusion coefficient $a(x,y,z)$ at $z=0.8$. The small errors signify the effectiveness of the Laplace-fPINNs method in identifying the three-dimensional diffusion coefficient.  The second line of Figure~\ref{fig_3} compares the exact solution $u(x,y,z,t)$ with the approximate solution $u_{NN}(x,y,z,t;\theta_{1}^{*})$ at $z=0.8$, $t=1$. The results indicate that the Laplace-fPINNs method can reconstruct the diffusion coefficient $a$ and predict the solution $u$ concurrently by incorporating additional measurement data $h(x,y,z,t)$.

\section{Conclusion}\label{sec4}
In this paper, a deep learning method for solving forward and inverse problems in subdiffusion is provided. Using the Laplace transform of time-fractional derivative, physics-informed neural networks, and numerical inverse Laplace transform, we propose a Laplace-fPINNs method that can infer a solution from the time-fractional diffusion equation (\ref{1.1}). The resulting method essentially avoids the use of auxiliary points, which were introduced in~\cite{Pang+Lu+Karniadakis-2019} to discretize the time-fractional Caputo derivative. A series of numerical results demonstrate the feasibility of this method for solving the forward problem of the subdiffusion. Additionally, we apply this method to solve an inverse problem that aims to determine the three-dimensional diffusion coefficient, which is challenging to solve using traditional inversion techniques due to the problem's high dimensionality. Finally, we remark that the idea of the Laplace-fPINNs method can be easily extended to the physics-informed operator learning methods for solving the subdiffusion
equation (\ref{1.1}).

\bigskip
\noindent{\bf Acknowledgments}
This work is sponsored by the National Key R\&D Program of China  Grant No. 2022YFA1008200 (Z. X.) and No. 2020YFA0712000 (Z. M.), the Shanghai Sailing Program (Z. X.), the Natural Science Foundation of Shanghai Grant No. 20ZR1429000  (Z. X.), the National Natural Science Foundation of China Grant No. 62002221 (Z. X.), the National Natural Science Foundation of China Grant No. 12101401 (Z. M.), the National Natural Science Foundation of China Grant No. 12031013 (Z. M.), Shanghai Municipal of Science and Technology Major Project No. 2021SHZDZX0102, and the HPC of School of Mathematical Sciences and the Student Innovation Center, and the Siyuan-1 cluster supported by the Center for High Performance Computing at Shanghai Jiao Tong University.
\medskip

\bibliographystyle{unsrt}
\bibliography{references}

\end{document}